\documentclass[letterpaper]{article}

\usepackage[T1]{fontenc}
\usepackage{multirow}
\usepackage{geometry}
\usepackage{setspace}

\usepackage[style = chem-acs, articletitle=true]{biblatex}
\usepackage{graphicx}
\usepackage{float}
\newfloat{scheme}{htbp}{los}
\floatname{scheme}{Scheme}
\floatname{chart}{Chart}
\newfloat{graph}{htbp}{loh}

\usepackage{chemformula} 
\usepackage[version = 4]{mhchem} 

\usepackage{amssymb}
\usepackage{amsmath}
\usepackage{mathrsfs}
\usepackage{algorithmicx}
\usepackage{algpseudocode}
\usepackage[utf8]{inputenc}
\usepackage{booktabs} 
\usepackage{siunitx}  
\allowdisplaybreaks
\usepackage[ruled, vlined, linesnumbered]{algorithm2e}
\usepackage{tikz}
\usetikzlibrary{decorations.pathreplacing, arrows.meta, positioning}
\usepackage{etoolbox}
\robustify\bfseries

\def \rd {\mathrm d}

\usepackage{authblk}
\author[1]{Sheng Chen}
\affil[1]{Research Center for Mathematics,  Beijing Normal University, Zhuhai, China.}

\author[2]{Sihong Shao}
\affil[2]{CAPT, LMAM and School of Mathematical Sciences, Peking University, Beijing, China.}

\author[3]{Shuai Wu*}
\affil[3]{School of Mathematical Sciences, Beijing Normal University, Beijing, China.}

\title{An Adaptive Log-Laguerre Spectral Method for the Radial Dirac Equation: Resolving Asymptotic Decay and Core Singularities in Atomic Calculations}
\date{*Email: ws@mail.bnu.edu.cn}

\begin{document}

\maketitle

\begin{abstract}
  The high-precision solution of the radial Dirac equation is fundamental to relativistic quantum chemistry, essential for reliable pseudopotential generation and all-electron electronic structure methods. Capturing both the non-polynomial singularities at the origin and the state-dependent asymptotic decay on semi-infinite domains presents a significant computational challenge. In this work, we propose the Adaptive Log-Laguerre Spectral Method (ALLSM), a novel coupled spectral-element solver that seamlessly integrates three advanced mathematical methodologies into a unified framework. Specifically, Generalized Log-Orthogonal Functions (GLOFs) are deployed in the near-core region to intrinsically approximate complex $r^s$ singular behaviors without requiring prior knowledge of the exact analytical exponent $s$. Concurrently, an adaptive Laguerre spectral method is employed to dynamically capture diverse exponential tails on $[0, \infty)$, avoiding artificial domain truncation. To structurally guarantee spectral purity across this bipartite basis, the framework rigorously incorporates the Inverse Dirac Operator Method (IDOM), effectively eliminating variational collapse and spurious states. Validated across diverse physical regimes, including Coulomb, finite-nucleus, and screened potentials, the proposed solver restores exponential convergence and consistently achieves relative accuracies of $10^{-10}$. This work provides a robust, pollution-free computational kernel for atomic structure calculations, establishing a highly reliable numerical standard for complex molecular simulations.
\end{abstract}

\section*{Keywords}
Adaptive Laguerre methods, Radial Dirac equation, High-precision computation, Galerkin method, Log-Orthogonal functions, Spectral pollution, Singular problems.

\section*{Abbreviations}

\noindent
\begin{tabular}{@{}ll}
	\textbf{RDEQ}  & Radial Dirac Equation \\
	\textbf{ALLSM} & Adaptive Log-Laguerre Spectral Method \\
	\textbf{LOFs}  & Log-Orthogonal Functions \\
	\textbf{GLOFs} & Generalized Log-Orthogonal Functions \\
	\textbf{LFs}   & Laguerre Functions \\
	\textbf{GLFs}  & Generalized Laguerre Functions \\
	\textbf{IDOM}  & Inverse Dirac Operator Method \\
\end{tabular}

\section{Introduction}
The Radial Dirac Equation (RDEQ) is the heart of high-fidelity atomic structure calculations, serving as the indispensable computational kernel for modern relativistic electronic structure methods \cite{wilson2006purgatorio}, including Density Functional Theory \cite{hohenberg1964inhomogeneous} and ab initio pseudopotential generation \cite{jonsson2007grasp2k, salvat1991accurate}. Traditionally, the RDEQ is tackled using either shooting methods \cite{vcertik2013dftatom, press2007numerical, salvat1991accurate} or matrix diagonalization via basis set expansions, utilizing B-splines or finite elements \cite{pask2005finite, grant2009b, fischer2009b, vcertik2024high}. However, achieving high precision within a short time frame remains computationally demanding, primarily due to the intricate mathematical behavior of the radial wavefunctions. Numerical schemes are primarily challenged by two multi-scale spatial features: the state-dependent asymptotic decay on the far-field domain, and the non-polynomial behavior near the origin—particularly the severe $r^s$ singularities characteristic of point-nucleus Coulomb potentials, or the steep near-core gradients in finite-nucleus models. Moreover, basis set approaches must handle the intrinsic risk of variational collapse \cite{schwarz1982mol, mark1982new, kutzelnigg1984basis, lewin2010spectral}. Failing to rigorously resolve these challenges, while maintaining spectral stability, compromises the numerical reliability of atomic simulations.

Fundamentally, addressing these spatial bottlenecks is a challenge of function approximation. At the far-field boundary, bound-state wavefunctions exhibit exponential decay over the semi-infinite physical domain $[0, \infty)$, with precise decay rates and effective supports that are strictly state-dependent. Modern finite element and B-spline methods effectively resolve these challenges through automated box-size convergence and adaptive mesh grading \cite{vcertik2024high}. However, capturing highly oscillatory near-field nodes alongside diffuse far-field tails across diverse eigenstates often necessitates dense localized grids or state-specific mesh optimizations. Spectral methods offer a compelling alternative, capable of achieving exponential convergence with a highly compact degree of freedom (DOF) footprint, provided the basis functions strictly match the underlying physical asymptotics. Yet, standard global spectral bases struggle with multi-scale decay and are structurally ill-equipped to resolve non-polynomial near-core behaviors without losing spectral accuracy.

To systematically overcome these interconnected bottlenecks while maintaining an extremely compact basis representation, we propose an adaptive spectral-element framework: the Adaptive Log-Laguerre Spectral Method (ALLSM). Because a structurally stable computational kernel is paramount, our framework natively incorporates the Inverse Dirac Operator Method (IDOM) \cite{hill1994solution} to rigorously circumvent variational collapse. Built upon this variationally stable foundation, we employ a bipartite domain decomposition that strategically integrates two specialized functional families:
\begin{itemize}
	\item \textbf{Adaptive Laguerre method}: To capture the state-dependent asymptotic decay on the semi-infinite domain without heuristic domain truncation, we employ an adaptive Laguerre method \cite{xia2021efficient}. By treating the frequency indicator of the basis functions as an optimization parameter, this approach automatically aligns the basis scaling with the physical wavefunction's far-field tail, enabling standard orthogonal families to efficiently represent diverse quantum states.
	\item \textbf{Log-Orthogonal Functions}: To resolve non-polynomial behaviors near the origin, we deploy GLOFs \cite{chen2022log} in the near-core element. Crucially, the GLOF basis provides a ``black-box'' capability to intrinsically approximate complex singular or steep gradient behaviors, such as the fractional $r^s$ dependence, without requiring any prior analytical extraction of the exponent $s$.
\end{itemize}

Although the individual mathematical components utilized in this study—GLOFs, adaptive frequency-indicator-based scaling, and IDOM—have been established in previous literature \cite{hill1994solution, xia2021efficient, chen2022log}, our work introduces a unique computational paradigm for the RDEQ through the following contributions: (i) Robust Framework Integration: We propose a unified, variationally stable spectral-element framework that seamlessly integrates the near-core non-polynomial treatment of GLOFs with the far-field adaptive Laguerre methodology via a mathematically rigorous $C^0$-continuous bridging construction. This approach resolves the multi-scale spatial challenges of the Dirac equation as a closed-loop system. (ii) Structural Pollution Elimination: We provide a rigorous topological analysis—rooted in the Levitin-Shargorodsky theorem—explaining how the inverse operator formulation reconfigures the essential spectrum to provide a pollution-free approximation in all tested cases, ensuring absolute spectral purity regardless of upper-lower DOF imbalances. (iii)High-Fidelity Benchmarking: Validated across diverse regimes, including point-Coulomb, Gaussian-nucleus, and screened interactions, this compact framework generates highly stable reference data for complex or mathematically intractable effective potentials, serving as a rigorous numerical standard for benchmarking relativistic algorithms.

The remainder of this paper is organized as follows. Section \ref{sec: 2} establishes the theoretical foundation, introducing the continuous RDEQ, the exact asymptotic behavior of the wavefunctions, and the approximation methodologies (GLOFs and adaptive Laguerre). Section \ref{sec: 3} provides a rigorous theoretical dissection of spectral stability, elucidating the topological origins of variational collapse and justifying the justifying the use of IDOM. Section \ref{sec: 4} outlines the practical computational implementation, detailing the globally continuous spectral basis and the Galerkin weak formulation. In Section \ref{sec: 5}, we comprehensively evaluate the proposed framework across a diverse array of chemical and physical potentials. Finally, Section \ref{sec: 6} summarizes our core methodological contributions and outlines trajectories for future applications.

\section{Theoretical Framework and Approximation Methodology}\label{sec: 2}
The theoretical foundation of our computational framework is the time-independent Dirac equation for a single particle in a central potential $V(r)$. Taking advantage of spherical symmetry, the four-component spinor eigenstates decouple, reducing the problem to the Radial Dirac Equation (RDEQ)
\begin{equation}\label{Dirac}
	\mathbf{D}_r \Psi = 
	\begin{pmatrix} 
		V(r) + c^2 & -c\left(\frac{\rd}{\rd r} - \frac{\kappa}{r}\right) \\[5pt]
		c\left(\frac{\rd}{\rd r} + \frac{\kappa}{r}\right) & V(r) - c^2 
	\end{pmatrix} 
	\begin{pmatrix} 
		F_{n\kappa}(r) \\[5pt]
		G_{n\kappa}(r) 
	\end{pmatrix} 
	= 
	E_{n\kappa} 
	\begin{pmatrix} 
		F_{n\kappa}(r) \\[5pt]
		G_{n\kappa}(r) 
	\end{pmatrix},
\end{equation}
where $F_{n\kappa}(r)$ and $G_{n\kappa}(r)$ are the large and small radial components, respectively. The energy level, $E_{n\kappa}$, is solely dependent on the principal quantum number $n$ and the total angular momentum quantum number $\kappa$. The normalization condition is given by 
\begin{equation*}
	\int_{0}^{\infty} (|F_{n\kappa}(r)|^2 + |G_{n\kappa}(r)|^2) ~\mathrm{d}r = 1.
\end{equation*}

\subsection{The Asymptotic Properties and Spatial Challenges}
To understand the fundamental problems in numerically solving the RDEQ, it is necessary to analyze the spatial behavior of the wavefunctions \cite{zabloudil2005electron}. For simplicity, we will omit the energy levels during derivation and keep the concise notation $F_{n\kappa}\to F,~G_{n\kappa}\to G,~E_{n\kappa}\to E.$ 

\subsubsection*{Asymptotic Properties Near the Origin}
Near the origin, for models exhibiting an infinite singularity, such as the classic point-nucleus Coulomb potential where $\lim_{r\to 0} V(r) \sim -Z/r$, asymptotic analysis dictates that both radial components behave as $r^s$ at the origin:
\begin{equation*}
	\lim_{r \to 0} F(r) \sim r^s, \quad \lim_{r \to 0} G(r) \sim r^s, \quad \text{where } s = \sqrt{\kappa^2 - (Z/c)^2}.
\end{equation*}
Because the relativistic angular momentum quantum number $\kappa$ is a non-zero integer, the exponent $s$ is strictly non-integer. This non-polynomial behavior significantly slows the convergence rates of classical polynomial-based spectral methods. While modern high-order finite element methods handle this effectively via localized geometric mesh refinement \cite{vcertik2024high}, designing a specialized global or semi-global basis to intrinsically capture such behavior within an extremely compact DOF remains a significant mathematical hurdle. 

Furthermore, for physically realistic finite-nucleus models, e.g., Gaussian distributions, the potential is regularized at the origin, eliminating the $r^s$ singularity. However, the wavefunction still exhibits steep, multi-scale gradients near the nuclear boundary. A robust computational framework must handle both the unknown fractional exponents of singular models and the steep near-core variations of regularized models efficiently without requiring a priori analytical extractions.

\subsubsection*{Asymptotic Properties for the Far Field}
\begin{figure}[htbp]
	\begin{center}
		\includegraphics[scale=0.45]{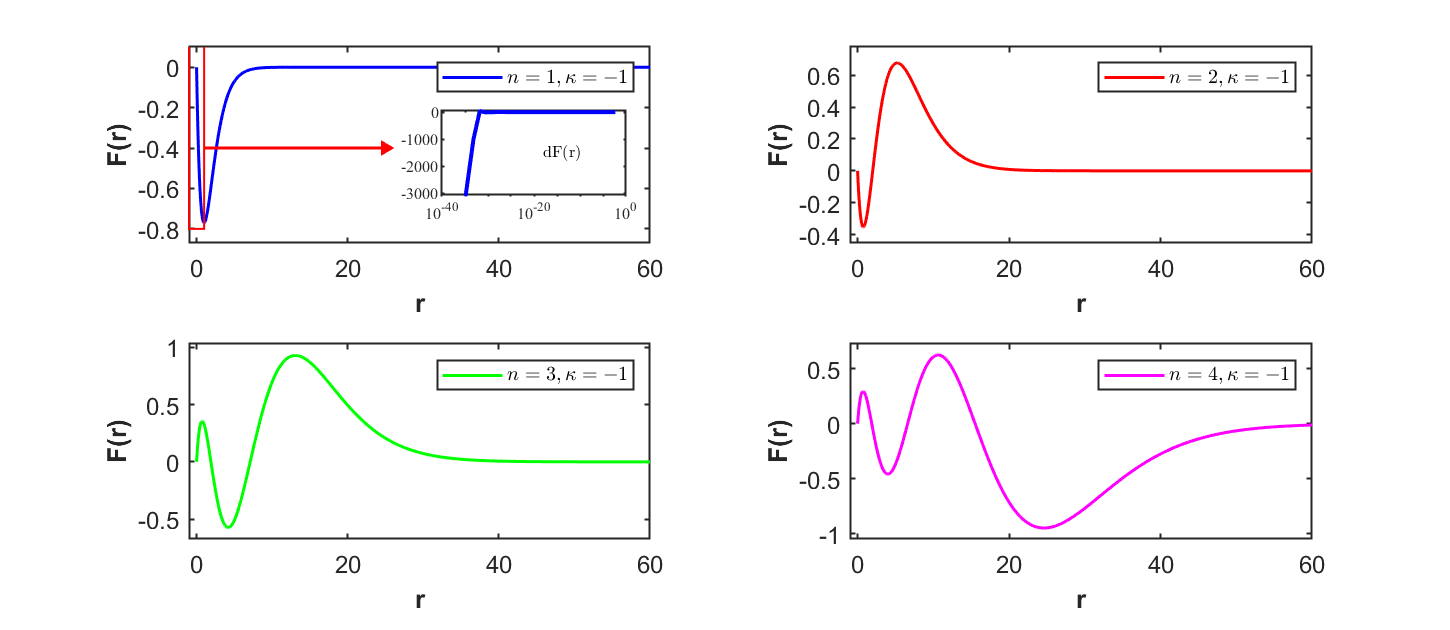}
	\end{center}
	\caption{{Comparison of relativistic large-component radial wavefunctions $F(r)$ and their derivatives for hydrogenic potentials, $n=1, \dots, 4, \kappa = -1$. (a) The top-left panel highlights the singular nature near the origin; its inset plots the derivative $\mathrm dF(r)$ on an ultra-fine $10^{-30}$ radial scale, revealing its rapid divergence to infinity. (b) The other panels illustrate the diverse decay rates in the far-field and oscillations in the near-field across different principal quantum numbers $n$.}}
	\label{fig: 10} 
	\begin{center}
		\includegraphics[scale=0.45]{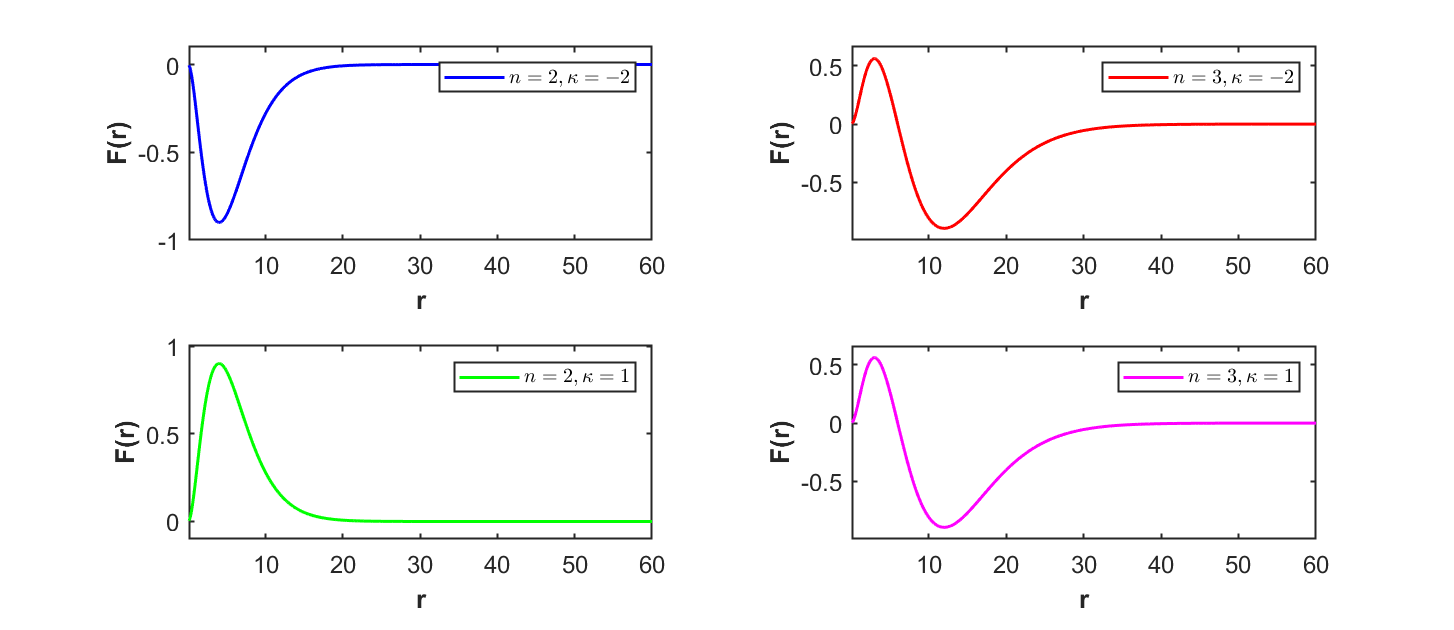}
	\end{center}
	\caption{{Comparison of relativistic large-component radial wavefunctions $F(r)$ for hydrogenic potentials with varying quantum numbers, $n=2, 3$ and $\kappa = -2, 1$. The individual panels illustrate the structural evolution of the wavefunctions across different relativistic states. Specifically, they showcase the distinct nodal structures and near-field oscillations determined by the principal quantum number $n$ and the relativistic angular momentum quantum number $\kappa$.}}
	\label{fig: 12} 
\end{figure}

Conversely, in the far-field limit where the potential vanishes, $\lim_{r \to \infty} V(r) = 0$, the RDEQ asymptotically decouples into a second-order form: 
\begin{equation*}
	\frac{\rd^2  F}{\rd r^2} \approx \frac{1}{c}(c^2-E)\frac{\rd G}{\rd r} = \frac{1}{c}(c^2-E) \Big[ ~\frac{1}{c}(E+c^2) F~ \Big] = \frac{c^4-E^2}{c^2}  F.
\end{equation*}
Governing the exponential decay of the bound states, we have 
\begin{equation*}
	\lim_{r \to \infty} F(r) \sim e^{-\lambda r}, \quad \text{where } \lambda = \sqrt{c^2 - (E/c)^2}.
\end{equation*}
Crucially, this decay rate $\lambda$ is strictly state-dependent. As illustrated in Figures \ref{fig: 10} and \ref{fig: 12}, the effective spatial support and near-field nodal structures of the wavefunctions vary drastically across different quantum states. Rigid, fixed-scale global basis sets are structurally inadequate to simultaneously capture drastic near-core fluctuations and diverse far-field tails efficiently. Therefore, a specialized basis capable of adaptively optimizing its spatial resolution is essential.

\subsection{Log-Orthogonal Functions for Non-Polynomial Singularities}
By the asymptotic analysis, the RDEQ solutions exhibit a non-polynomial fractional power behavior, $r^s$, in the vicinity of the origin. To intrinsically capture this singularity without resorting to cumbersome geometric mesh refinements or explicit analytical extractions, we deploy GLOFs \cite{chen2022log, chen2020spectrally} in the near-origin region. 

The standard LOFs are defined via a logarithmic mapping to classical Laguerre polynomials, 
\begin{equation*}\label{SLF}
	\mathcal{S}_n^\gamma(x) = \mathscr{L}_n(-(\gamma+1)\log{x}), \quad x \in (0,1),
\end{equation*}
where $\mathscr{L}_n$ denotes the Laguerre polynomial of degree $n$, and the parameter $\gamma > -1$ regulates the resolution of distinct singular behaviors. While standard LOFs are theoretically important, their direct application in Galerkin formulations often suffers from numerical ill-conditioning. To guarantee structural stability and ensure the basis functions strictly reside within the $L^2$ space, we introduce a practically optimized variant known as the GLOFs\cite{chen2022log}, 
\begin{equation*}
	\mathcal S_n^{(\gamma,\mu)}(x) = x^{(\gamma-\mu)/2}\mathcal S_n^{\gamma}(x),\quad n\geq 0.
\end{equation*}
For relativistic atomic structure calculations, setting $\mu = 2$ and $\gamma = 4$ proves to be a proper choice\cite{chen2025efficient}. This yields the explicitly well-conditioned GLOF basis,
\begin{equation*}\label{GLOFs}
	\mathcal S_n^{(4,2)}(x) = x\mathscr{L}_n(-5\log{x}).
\end{equation*}

Furthermore, constructing a homogeneous basis satisfying strict zero boundary conditions at both $x=0$ and $x=1$ is elegantly constructed via the linear combination,
\begin{equation*}
	\mathcal S_{n+1}^{(4,2)}(x)-\mathcal S_n^{(4,2)}(x) = x\big(\mathscr{L}_{n+1}(-5\log{x})-\mathscr{L}_n(-5\log{x})\big).
\end{equation*}
The computational advantage of the GLOFs is their automatic capability to approximate singular terms without any prior specification of the exact physical exponents. As established by the approximation theory of LOFs \cite{chen2022log}, if a generic multi-scale wavefunction is modeled as $u(x) = \sum_{i=1}^K x^{s_i}$ with completely unknown exponents $s_i > 0$, the convergence rate under the $L^2$ projection $\Pi^{4,2}_N$ is strictly governed by the ratio,
\begin{equation*}
	R_{4,2} = \max_{1\leq i\leq K}\left| \frac{2s_i + 2 - 4}{2s_i + 2 + 2 + 4} \right| = \max_{1\leq i\leq K}\left| \frac{s_i -1}{s_i +4} \right| < 1.
\end{equation*}
The corresponding approximation error decays exponentially as,
\begin{equation*}
	|u - \Pi_N^{4,2} u|{x^2} \leq K\sqrt{\frac{25}{4+s}} N^{\frac{1}{2}} (R_{4,2})^N, \quad \text{for}\quad  N\geq 1.
\end{equation*}
Consequently, this tailored GLOF basis effectively resolves the severe $r^s$ non-polynomial asymptotics natively. More importantly, this establishes a robust mathematical framework that can be seamlessly extended to more complex differential equations in computational chemistry, where the precise analytical nature of the near-core singularity may be totally intractable.

\subsection{Asymptotic Capture via the Adaptive Laguerre Method}
\begin{figure}[ht]
	\centering
	\begin{minipage}{0.48\linewidth}
		\centering
		\includegraphics[width=\linewidth]{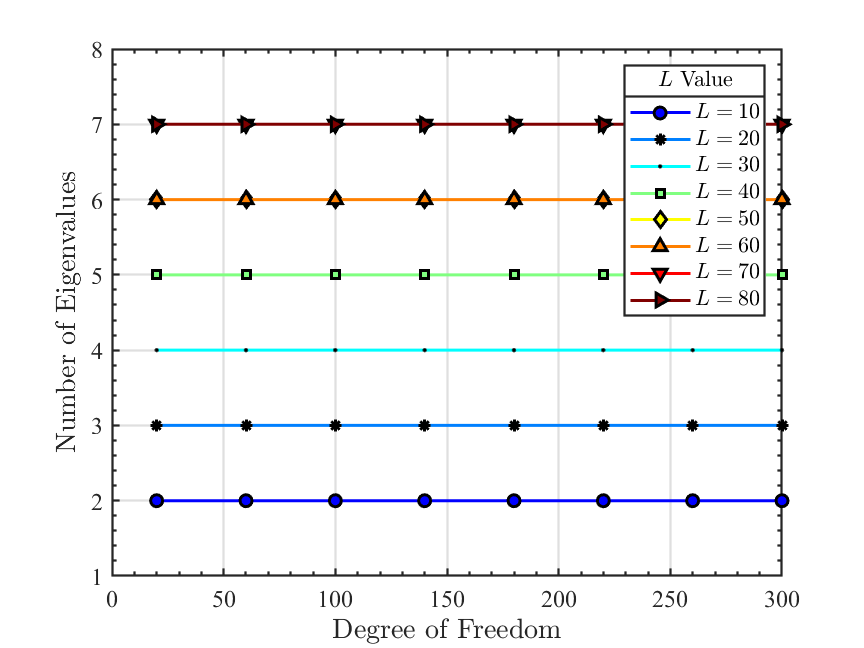}
	\end{minipage}
	\hfill
	\begin{minipage}{0.48\linewidth}
		\centering
		\includegraphics[width=\linewidth]{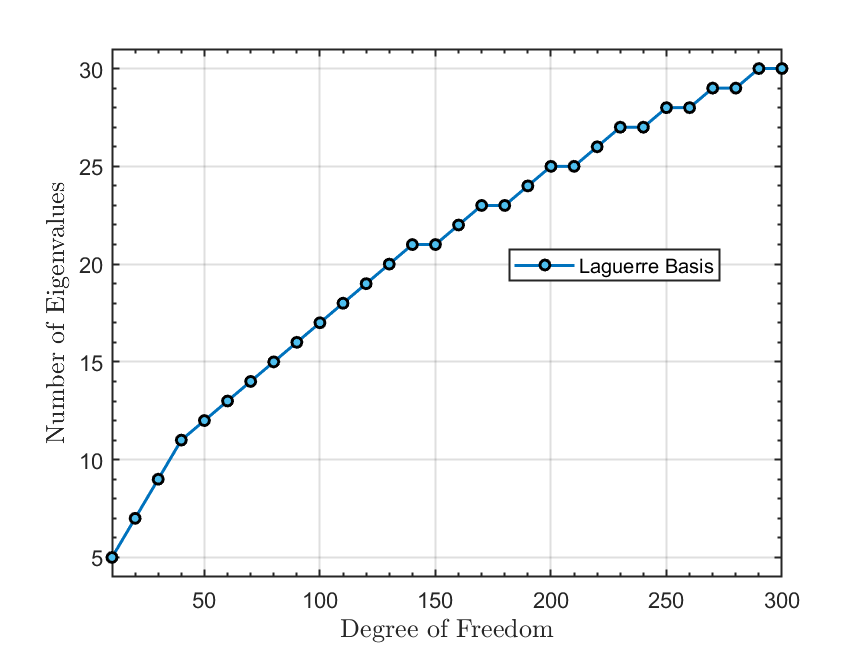}
	\end{minipage}
	\caption{The number of converged physical eigenvalues\quad \textbf{Left}: Legendre basis. \quad \textbf{Right}: Laguerre basis.}
	\label{fig: 131}
	\begin{minipage}{1.0\linewidth}
		\centering
		\includegraphics[width=\linewidth]{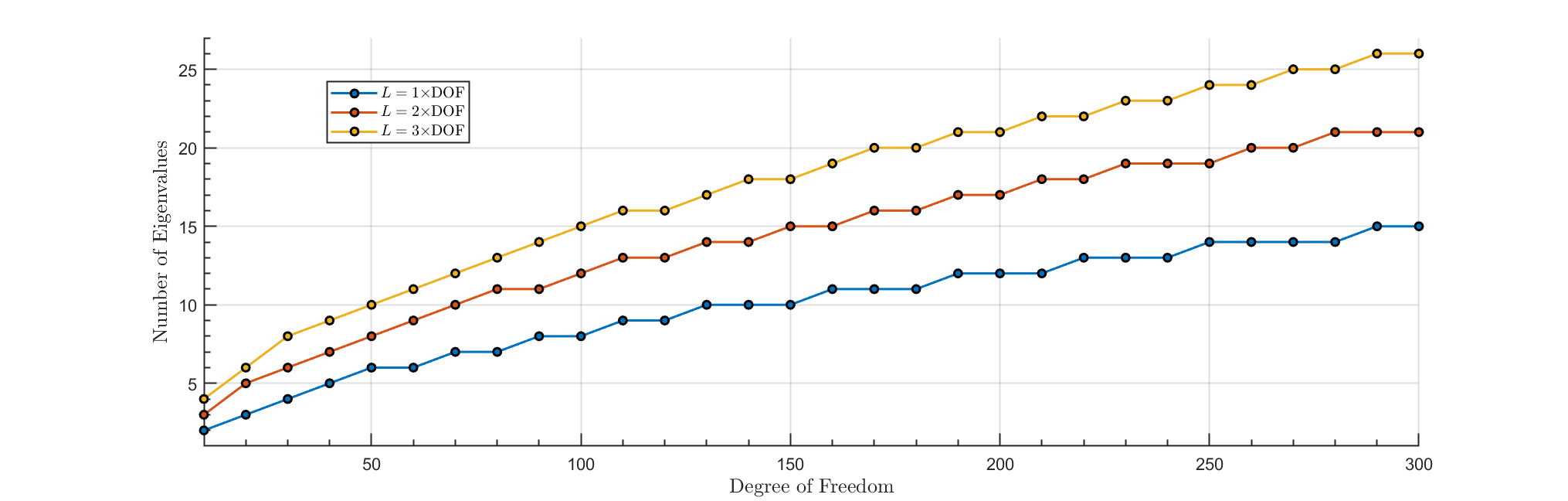}
	\end{minipage}
	\caption{The number of eigenvalues with degree of freedom for fixed ratio.}
	\label{fig: 1311}
\end{figure}

When employing a truncated basis defined on a finite interval $(0, L)$, such as Legendre polynomials, the quantum system is effectively confined within an artificial bounding box. While modern numerical algorithms robustly handle this by systematically extending the domain $L$ and proportionally increasing the basis size $N$ to maintain resolution density, capturing highly diffuse, higher-energy states introduces an inherent computational overhead. As illustrated in Figures \ref{fig: 131} and \ref{fig: 1311}, the number of converged physical eigenvalues is intricately linked to the truncation boundary. To perfectly envelop states with broader spatial support, both the domain size and the degrees of freedom must be continuously inflated.

By contrast, for quantum mechanical problems defined on a semi-unbounded domain $[0, \infty)$, Laguerre basis functions emerge as a mathematically natural choice that bypasses the need for spatial confinement. Because the Laguerre family intrinsically possesses an exponential tail, it naturally accommodates the expansive spatial domains of relativistic wavefunctions. Consequently, increasing the basis size $N$ directly translates into a higher-fidelity resolution of the physical spectrum without the necessity of adjusting an artificial boundary. 

However, a fixed-scale Laguerre basis is ill-suited for the multi-scale, state-dependent decay discussed previously. Capturing both highly oscillatory near-field nodes and diffuse far-field tails simultaneously still requires a prohibitive number of degrees of freedom. To fully exploit the Laguerre tail and restore extreme exponential convergence, we introduce a dynamically optimized scaling parameter $\beta$. This yields the Generalized Laguerre Functions (GLFs) $\widehat{\mathscr{L}}_n^\beta(r)$, which effectively scale the coordinate system to perfectly match the physical support of the target eigenstate.

Inspired by adaptive spectral methods \cite{xia2021efficient}, we treat the scaling parameter $\beta$ as an optimization variable governed by a \emph{frequency indicator}, $\mathcal{F}$. This indicator quantifies the relative energy content residing in the high-frequency tail of the computed eigenfunction $U_N^{\beta}(r) = \sum_{l=0}^N c_l^{\beta} \widehat{\mathscr{L}}_l^\beta(r)$. A minimized $\mathcal{F}$ implies rapid decay of the spectral coefficients, guaranteeing high-fidelity asymptotic representation. The indicator is rigorously defined as:
\begin{equation*}
	\mathcal{F}(U_N^{\beta}) = \left( \frac{\sum_{l=N-M+1}^{N} \gamma_l^{\beta} (c_l^{\beta})^2}{\sum_{l=0}^{N} \gamma_l^{\beta} (c_l^{\beta})^2} \right)^{\frac{1}{2}},
\end{equation*}
where $c_l^{\beta}$ are the spectral coefficients, $\gamma_l^{\beta} = \langle\widehat{\mathscr{L}}_n^\beta(r)|\widehat{\mathscr{L}}_n^\beta(r)\rangle$ are the weighting factors, and $M = \big\lfloor\frac{N}{3}\big\rfloor$ isolates the upper third of the spectrum for dealiasing \cite{hou2007computing}. To establish an automated stopping criterion, we define a reference threshold $f_0 = \mathcal{F}(U_M)$ derived from the well-behaved model function $U_M(r) = \sin(r)\exp(-r)$ \cite{shen2009some}. Because $U_M$ typically achieves machine precision $10^{-16}$ at $N=80$, its frequency indicator $f_0$ serves as a rigorous benchmark.

The adaptive scheme systematically optimizes the scaling parameter $\beta$ to minimize $\mathcal{F}(U_N^\beta)$, as detailed in Algorithm \ref{alg: freq-scaling}. This fully automated workflow allows the computational kernel to seamlessly adapt to diverse quantum states, effectively functioning as a black-box optimizer without requiring manual tuning or prior physical intuition regarding the specific decay rate.
\begin{algorithm}[h]
	\SetKwInOut{Input}{Input}
	\SetKwInOut{Output}{Output}
	\caption{Adaptive Frequency-Based Scaling for Eigenvalue Problems}
	\label{alg: freq-scaling}
	\Input{Truncation degree $N$, expansion factor $\nu > 1$, lower bound $\underline{\beta}$, upper bound $\overline{\beta}$, desired eigenmode index $i$.}
	\Output{Optimal parameter $\beta$, target characteristic value $E_i$, and the $i$th eigenfunction $U_{N}^{\beta}(x)$.}
	
	\BlankLine
	\KwSty{Define} $\text{COMPUTE\_EIGENFUNCTION}(i, \beta_k)$ as solving the GEP for the $i$-th eigenpair using scaling $\beta_k$.\\
	\KwSty{Define} $\text{FREQUENCY\_INDICATOR}(U)$ as computing $\mathcal{F}(U)$.\\
	\KwSty{Initialize} $\nu > 1,~ \beta=\beta_0,~ \underline{\beta}=10^{-4},~ \overline{\beta}=10^5,~ U_M \leftarrow \sin(x)\exp(-x)$.
	
	$U_{N}^{\beta} \leftarrow \text{COMPUTE\_EIGENFUNCTION}(i, \beta)$ \hfill // Compute the initial eigenfunction approximation\;
	$f \leftarrow \text{FREQUENCY\_INDICATOR}\big(U_{N}^{\beta}\big) $\;
	$f_0 \leftarrow \text{FREQUENCY\_INDICATOR}\big(U_M\big) $\;
	
	\If{$f > \nu f_0$}
	{
		$\beta_{\text{down}} \leftarrow \frac{1}{\nu}\beta$\;
		$U_{\text{down}} \leftarrow \text{COMPUTE\_EIGENFUNCTION}(i, \beta_{\text{down}})$\;
		$f_{\text{down}} \leftarrow \text{FREQUENCY\_INDICATOR}\big(U_{\text{down}}\big)$\;
		
		$\beta_{\text{up}} \leftarrow \nu\beta$\;
		$U_{\text{up}} \leftarrow \text{COMPUTE\_EIGENFUNCTION}(i, \beta_{\text{up}})$\;
		$f_{\text{up}} \leftarrow \text{FREQUENCY\_INDICATOR}\big(U_{\text{up}}\big)$\;
		
		$q \leftarrow \nu \cdot \mathbb{I}(f_{\text{up}} > f_{\text{down}}) + \frac{1}{\nu} \cdot \mathbb{I}(f_{\text{up}} \leq f_{\text{down}})$ \hfill // Determine the scaling direction $q$ ($\mathbb{I}$ is the indicator function)\;
		
		$\tilde{\beta} \leftarrow q\beta$\;
		$U_{N}^{\tilde{\beta}} \leftarrow \text{COMPUTE\_EIGENFUNCTION}(i, \tilde{\beta})$\;
		$\tilde{f} \leftarrow \text{FREQUENCY\_INDICATOR}\big(U_{N}^{\tilde{\beta}}\big)$\;
		
		\While{$\tilde{f} \le f$ \textbf{and} $(\tilde{\beta} \ge \underline{\beta}$ \textbf{or} $\tilde{\beta} \le \overline{\beta})$}
		{
			$\beta \leftarrow \tilde{\beta}$\;
			$U_{N}^{\beta} \leftarrow U_{N}^{\tilde{\beta}}$\;
			$f \leftarrow \tilde{f}$\;
			$\tilde{\beta} \leftarrow q\beta$\;
			$U_{N}^{\tilde{\beta}} \leftarrow \text{COMPUTE\_EIGENFUNCTION}(i, \tilde{\beta})$\;
			$\tilde{f} \leftarrow \text{FREQUENCY\_INDICATOR}(U_{N}^{\tilde{\beta}})$\;
		}
	}
\end{algorithm}

\section{Topological Origins of Spectral Stability}\label{sec: 3}
Having established the spatial approximation methodologies, we now dissect the fundamental numerical mechanics governing spectral stability. Solving the discrete Dirac equation is notoriously fraught with instabilities, most notably the phenomenon of variational collapse. In this section, we investigate the relationship between degree of freedom mismatches and the emergence of spurious states, providing a rigorous theoretical justification for the stabilizing supremacy of the inverse operator formulation.

\subsection{Variational Collapse and DOF Mismatch}
Variational collapse has remained a notorious pathology in relativistic quantum chemistry since its early diagnosis \cite{schwarz1982mol, schwarz1982two}. Historically, this emergence of spurious states was attributed to an intuitive physical cause: the Dirac Hamiltonian $\hat{\mathbf{D}}$ is not bounded from below, ostensibly causing variational procedures to ``fall'' into the negative-energy Dirac sea \cite{schwarz1982mol, mark1982new, kutzelnigg1984basis}. Conventional remedies, such as the kinetic balance condition, were developed primarily to enforce the correct non-relativistic Schrödinger limit and artificially stave off this collapse. However, this pragmatic intuition obscures the rigorous mathematical root of the instability.

The true mechanism of variational collapse is strictly topological \cite{lewin2010spectral}. Because the essential spectrum of $\hat{\mathbf{D}}$ is fundamentally non-connected, $(-\infty, -c^2] \cup [c^2, \infty)$, the physical bound-state mass gap $(-c^2, c^2)$ is left completely topologically exposed. In standard finite basis representations, it is precisely this structural vulnerability that allows spurious eigenvalues to freely invade the gap and contaminate the physical spectrum. Realistic molecular potentials that vanish at infinity inevitably resurrect the non-connected Dirac sea, dictating that numerical stability must rely on a mathematically rigorous discretization scheme rather than empirical physical constraints.

A critical vector for this spectral contamination is the imbalance in the degrees of freedom between the upper spinor components $F$ (with $N$ basis functions) and the lower spinor components $G$ (with $K$ basis functions). Our numerical investigations reveal a highly directional pollution mechanism governed by this mismatch. When the upper DOF exceeds the lower ($N > K$), spurious solutions predominantly descend from the positive continuum, clustering densely within the mass gap and obscuring the physical bound states, as evidenced in Table \ref{tab:results_data1}. Conversely, an excess of lower-component DOF ($K > N$) forcefully drives spurious states upward. Remarkably, in the perfectly balanced $N = K$ configuration, no spurious states are observed across our tested potentials using our compact spectral representation. This is a direct testament to the extreme efficiency of ALLSM, in which the required basis size is simply too compact to provide the topological flexibility necessary for variational collapse to manifest.

\begin{table}[ht]
	\centering
	\caption{Spectral pollution under degree of freedom mismatch, where $N=200$ and varying $K$. Bold values indicate spurious states. Note that as $K$ approaches $N$, the spurious states are progressively repelled from the physically relevant bound states, until they completely vanish in the perfectly balanced $N=K$ configuration.}
	\label{tab:results_data1}
	\small
	\sisetup{
		round-mode=places,
		round-precision=5,
		table-format=-1.5,
		detect-weight=true 
	}
	\begin{tabular}{l *{6}{S}}
		\toprule
		& {$K=100$} & {$K=120$} & {$K=140$} & {$K=160$} & {$K=180$} & {$K=200$} \\
		& {($N=200$)} & {($N=200$)} & {($N=200$)} & {($N=200$)} & {($N=200$)} & {($N=200$)} \\
		\midrule
		$n=1$ & \textbf{-2.65078} & \textbf{-1.55781} & \textbf{-0.82504} & -0.50001 & -0.50001 & {-0.50001} \\
		$n=2$ & \textbf{-1.52727} & \textbf{-0.93142} & \textbf{-0.50719} & \textbf{-0.36255} & -0.12500 & -0.12500 \\
		$n=3$ & \textbf{-0.55156} & -0.50001 & -0.50001 & \textbf{-0.22105} & \textbf{-0.10292} & -0.05556 \\
		$n=4$ & -0.50001 & \textbf{-0.33104} & \textbf{-0.17879} & -0.12500 & -0.05556 & -0.03125 \\
		$n=5$ & \textbf{-0.41387} & \textbf{-0.24924} & \textbf{-0.13329} & \textbf{-0.07899} & \textbf{-0.05452} & -0.02000 \\
		$n=6$ & \textbf{-0.23199} & \textbf{-0.13968} & -0.12500 & \textbf{-0.05678} & -0.03125 & -0.01389 \\
		$n=7$ & \textbf{-0.18866} & -0.12500 & \textbf{-0.07534} & -0.05556 & \textbf{-0.02094} & -0.01020 \\
		$n=8$ & \textbf{-0.12690} & \textbf{-0.11301} & \textbf{-0.06010} & \textbf{-0.03290} & -0.02000 & -0.00781 \\
		$n=9$ & -0.12500 & \textbf{-0.07637} & -0.05556 & -0.03125 & -0.01389 & -0.00617 \\
		$n=10$ & \textbf{-0.10741} & \textbf{-0.06421} & \textbf{-0.04108} & \textbf{-0.02555} & \textbf{-0.01389} & -0.00500 \\
		\bottomrule
	\end{tabular}
\end{table}

\subsection{Operator Reformulation: IDOM vs. SDOM}
Ultimately, we must mathematically demystify the robust stability guaranteed by operator transformation techniques. A pervasive misconception is that such methods prevent collapse merely by bounding the energy from below. The true mechanism of spectral stabilization is distinctly topological, governed by the Levitin-Shargorodsky theorem, which establishes that spectral pollution in Galerkin methods can manifest anywhere within the convex hull of the essential spectrum $\sigma_{ess}(\hat{\mathbf{D}})$. Because the original Dirac operator features essential spectra on both sides of the real axis, its convex hull spans the entire real line, rendering the entire mass gap a fertile breeding ground for spurious roots.

The effectiveness of the Inverse Dirac Operator Method (IDOM) and the Squared Dirac Operator Method (SDOM) lies entirely in their ability to radically reconfigure this essential topology. As shown in Figure \ref{fig:spectrum_mapping}, SDOM maps the dual spectra onto a single semi-infinite interval $[c^4, \infty)$, whereas IDOM compresses the unbounded spectra into a compact, connected interval $[-1/c^2, 1/c^2]$. Both transformations systematically collapse the interior of the essential spectrum's convex hull, mathematically destroying the topological gap required for pollution to form.

\begin{figure}[htbp]
	\centering
	\begin{tikzpicture}[>=stealth, scale=0.8, every node/.style={transform shape}]
		
		\tikzset{
			ess/.style={line width=4.5pt, blue!70!black},
			discrete/.style={circle, fill=red, inner sep=2pt},
			spurious/.style={circle, fill=gray!60, inner sep=2pt},
			hullbrace/.style={decorate, decoration={brace, amplitude=8pt, mirror}, thick, gray!80}
		}
		
		\begin{scope}[shift={(0,6)}]
			\node[font=\bfseries, align=center] at (0, 1) {Original Dirac Operator $\hat{\mathbf{D}}$};
			\draw[->, thick] (-5.5, 0) -- (5.5, 0) node[right, font=\large] {$E$};
			
			\draw[ess, <-] (-5, 0) -- (-2.5, 0) node[midway, above=3pt] {$\sigma_{ess}$};
			\draw[ess, ->] (2.5, 0) -- (5, 0) node[midway, above=3pt] {$\sigma_{ess}$};
			
			\draw[thick] (-2.5, 0.2) -- (-2.5, -0.2) node[below] {$-c^2$};
			\draw[thick] (2.5, 0.2) -- (2.5, -0.2) node[below] {$c^2$};
			
			\node[discrete, label=above:{\small Physical}] at (1.5, 0) {};
			\node[discrete] at (1.9, 0) {};
			\node[discrete] at (2.2, 0) {};
			
			\node[spurious, label=above:{\small Spurious}] at (-1.5, 0) {};
			\node[spurious] at (-0.8, 0) {};
			\node[spurious] at (0, 0) {};
			\node[spurious] at (0.6, 0) {};
			
			\draw[hullbrace] (5, -1) -- (-5, -1) 
			node[midway, below=5pt, text=black, align=center, font=\large] 
			{Convex Hull spans $\mathbb{R}$ \\ \textcolor{red}{\textbf{Mass Gap is Topologically Polluted}}};
		\end{scope}
		
		\begin{scope}[shift={(-5,2)}]
			\node[font=\bfseries, align=center] at (0, 1) {IDOM: $\hat{\mathbf{D}}^{-1}$};
			\draw[->, thick] (-3.5, 0) -- (4.5, 0) node[right, font=\large] {$\lambda$};
			
			\draw[ess] (-1.5, 0) -- (1.5, 0) node[midway, above=3pt] {$\sigma_{ess}$};
			
			\draw[thick] (-1.5, 0.2) -- (-1.5, -0.2) node[below] {$-1/c^2$};
			\draw[thick] (1.5, 0.2) -- (1.5, -0.2) node[below] {$1/c^2$};
			
			\node[discrete] at (2.5, 0) {};
			\node[discrete] at (3.1, 0) {};
			\node[discrete] at (3.8, 0) {};
			
			\draw[hullbrace] (1.5, -1.1) -- (-1.5, -1.1) 
			node[midway, below=5pt, text=black, align=center, font=\normalsize] 
			{Compact Convex Hull \\ \textcolor{green!50!black}{\textbf{Bound States Isolated}}};
		\end{scope}
		
		\begin{scope}[shift={(2.5,2)}]
			\node[font=\bfseries, align=center] at (2.5, 1) {SDOM: $\hat{\mathbf{D}}^2$};
			\draw[->, thick] (-1, 0) -- (5.5, 0) node[right, font=\large] {$\lambda$};
			
			\draw[ess, ->] (2.5, 0) -- (5, 0) node[midway, above=3pt] {$\sigma_{ess}$};
			
			\draw[thick] (2.5, 0.2) -- (2.5, -0.2) node[below] {$c^4$};
			
			\node[discrete] at (0.5, 0) {};
			\node[discrete] at (1.2, 0) {};
			\node[discrete] at (1.8, 0) {};
			
			\draw[hullbrace] (5, -1.1) -- (2.5, -1.1) 
			node[midway, below=5pt, text=black, align=center, font=\normalsize] 
			{Semi-infinite Hull \\ \textcolor{green!50!black}{\textbf{Bound States Isolated}}};
		\end{scope}
		
		\draw[->, thick, dashed, gray!80, bend right=15] 
		(-2.5, 5) to node[left=8pt, text=black, font=\large] {$E \mapsto 1/E$} (-5, 3.5);
		
		\draw[->, thick, dashed, gray!80, bend left=15] 
		(2.5, 5) to node[right=8pt, text=black, font=\large] {$E \mapsto E^2$} (5, 3.5);
		
	\end{tikzpicture}
	\caption{Topological stabilization mechanism of the Dirac spectrum via operator reformulations. \quad \textbf{Top}: The original Dirac operator $\hat{\mathbf{D}}$ features a convex hull spanning the entire real line, creating a ``breeding ground'' for spurious states within the mass gap. \quad \textbf{Bottom}: IDOM and SDOM reconfigure the essential spectrum $\sigma_{ess}$ into compact or semi-infinite intervals. By systematically collapsing the convex hull of the essential spectrum, these transformations destroy the topological conditions required for spectral pollution.}
	\label{fig:spectrum_mapping}
\end{figure}

To empirically demonstrate the superiority of IDOM over SDOM, we evaluate the relativistic harmonic oscillator ($V(r) = \frac{1}{2}r^2$). While SDOM eliminates the negative continuum, its non-injective mapping $E \mapsto E^2$ triggers \textit{spectral folding}. The true physical positive-energy states become densely interleaved with folded spurious states from the negative Dirac sea (Table \ref{table: kap=-11}), rendering automated algorithmic detection virtually impossible. In stark contrast, IDOM ($E \mapsto 1/E$) topologically isolates the target ground state as the absolute maximum positive eigenvalue, cleanly segregating the negative energy sea. Consequently, IDOM directly yields a pristine, reliably ordered spectrum without any spurious interleaving.

\begin{table}[!ht]
	\caption{Comparison of eigenvalues for the $\kappa=-1$ state in a Harmonic Oscillator potential. Bold values denote spurious states, while IDOM provides a clean, reliable benchmark for the physical states.}
	\label{table: kap=-11}
	\small 
	\centering	
	\begin{tabular}{l *{4}{S[table-format=2.8]}}
		\toprule
		& {IDOM} & {SDOM($N=30$)} & {SDOM($N=50$)} & {SDOM($N=70$)} \\ 
		\midrule
		$n=0$ & {2.49997504} & \textbf{0.06352360} & \textbf{1.01395279} & \textbf{1.34052631} \\  
		$n=1$ & {4.49983527} & 2.49998838 & 2.49997503 & 2.49997503 \\  
		$n=2$ & {6.49961565} & \textbf{3.26825736} & \textbf{2.82962861} & \textbf{3.01607031} \\  
		$n=3$ & {8.49931620} & 4.50045282 & 4.49983527 & 4.49983527 \\  
		$n=4$ & {10.49893692} & \textbf{6.51805642} & \textbf{4.85334683} & \textbf{4.85279944} \\  
		$n=5$ & {12.49847782} & 8.57215540 & 6.49961565 & 6.49961565 \\  
		$n=6$ & {14.49793893} & \textbf{8.58823544} & \textbf{7.12972519} & \textbf{6.87374186} \\  
		$n=7$ & {16.49732025} & \textbf{10.97222189} & 8.49931619 & 8.49931619 \\  
		$n=8$ & {18.49662179} & \textbf{14.19012901} & \textbf{9.68195672} & \textbf{9.08827958} \\  
		$n=9$ & {20.49584357} & \textbf{17.41947110} & 10.49893691 & 10.49893691 \\ 
		\bottomrule
	\end{tabular}
\end{table}

\section{Computational Details}\label{sec: 4}
In this section, we detail the practical implementation of the proposed adaptive spectral framework. To translate our theoretical strategies into a robust computational kernel, we proceed in three sequential steps: (1) constructing a unified globally continuous basis that seamlessly connects GLOFs with adaptively scaled Laguerre functions; (2) implementing the IDOM Galerkin formulation to structurally eliminate variational collapse; and (3) addressing the numerical conditioning and algebraic implementation of the resulting discrete system. 

\subsection{Bipartite Basis Construction}
To translate the spatial strategies formulated in Section \ref{sec: 2} into a computable discrete space, we employ a spectral-element domain decomposition method. The semi-unbounded physical space is partitioned into a singular near-core region $I_1 = (0,L)$ and an asymptotic decay region $I_2 = (L,\infty)$, as illustrated in Figure \ref{basis}. This bipartite decomposition allows us to deploy the GLOFs and the GLFs exactly where their respective approximation properties are physically required. 

Unlike empirical domain truncation strategies where the boundary determines the resolution of the far-field tail, the partition radius $L$ merely demarcates the near-core non-polynomial zone. As will be demonstrated in Section \ref{sec: 5}, the framework is highly robust to the choice of $L$, with $L=0.01$ serving as a consistently stable threshold across varied potentials, negating the need for heuristic parameter tuning.

\begin{figure}[htbp]
	\centering
	\begin{tikzpicture}[>=latex]
		\draw[->, thick] (0,0) -- (10.5,0) node[right] {$r$};
		
		\draw[thick] (0,0.1) -- (0,-0.1);
		\node[below, text height=2ex] at (0,-0.1) {$0$};
		
		\draw[thick, red] (3,0.1) -- (3,-0.1);
		\node[below, red, text height=2ex] at (3,-0.1) {$L$};
		
		\node[below, text height=2ex] at (10,-0.05) {$\infty$};
		
		\draw[<->, shorten <=2pt, shorten >=2pt, blue] (0,0.5) -- (3,0.5);
		\node[above, blue, font=\large] at (1.5,0.5) {$\boldsymbol{I_1 = (0, L)}$};
		
		\draw[->, shorten <=2pt, thick, green!60!black] (3,0.5) -- (10,0.5);
		\node[above, green!60!black, font=\large] at (6.5,0.5) {$\boldsymbol{I_2 = (L, \infty)}$};
		
		\node[text width=4.5cm, align=center, anchor=north] at (1.5,-0.7) {
			\textbf{Singular near-core region} \\
			\textcolor{blue}{Non-polynomial singularity} \\
			\vspace{1ex}
			\fbox{\begin{minipage}{4cm} \centering \textbf{\textcolor{blue}{Log-Orthogonal Functions}} \end{minipage}}
		};
		
		\node[text width=6cm, align=center, anchor=north] at (7.5,-0.7) {
			\textbf{Exponential decay region} \\
			\textcolor{green!60!black}{Far-field tail} \\
			\vspace{1ex}
			\fbox{\begin{minipage}{4cm} \centering \textbf{\textcolor{green!60!black}{Laguerre Functions}} \end{minipage}}
		};
	\end{tikzpicture}
	\caption{Schematic of the bipartite domain decomposition on the semi-unbounded interval. The physical space is partitioned into a singular near-core region $I_1$ and an asymptotic decay region $I_2$, where GLOFs and GLFs are respectively deployed.}
	\label{basis}
\end{figure}

To guarantee that the resulting wavefunctions are physically admissible—specifically, globally $C^0$-continuous across the entire domain—the discrete basis must be carefully stitched at the element interface $r = L$. We construct two sets of local homogeneous basis functions that strictly vanish at the interface $r=L$:
$$
\begin{aligned}
	\phi_{n}^1(r) &= 
	\begin{cases}
		r/L\big(\mathcal S^{(4,2)}_{1}(r/L)-\mathcal S^{(4,2)}_{{n}-1}(r/L)\big),\\[5pt]
		0,
	\end{cases} 
	&&
	\begin{aligned}
		r\in I_1;\\[5pt]
		r\in I_2.
	\end{aligned}
	\label{eq:basis_core} \\[5pt]
	\phi^2_{n}(r) &= 
	\begin{cases}
		0,\\[5pt]
		\widehat{\mathscr L}_{n}^\beta(r-L)-\widehat{\mathscr L}_{n-1}^\beta(r-L),
	\end{cases} 
	&&
	\begin{aligned}
		r\in I_1;\\[5pt]
		r\in I_2.
	\end{aligned}
	\label{eq:basis_tail}
\end{aligned}
$$
Here, the scaling parameter $\beta$ in Equation~\eqref{eq:basis_tail} is adaptively determined via the frequency indicator scheme (Algorithm \ref{alg: freq-scaling}) to automatically match the state-dependent exponential tail.

Because both local basis sets enforced with homogeneous Dirichlet conditions at $r=L$ intrinsically cannot represent a wavefunction possessing a non-zero amplitude at this boundary, we introduce a non-homogeneous ``bridging'' function $\phi^3(r)$. This function smoothly connects a linear interior profile to the fundamental Laguerre decay, ensuring strict $C^0$ continuity:
\begin{equation}\label{eq:hat}
	\phi^3(r) =
	\begin{cases}
		r/L, & r\in I_1,\\[3pt]
		e^{-\beta(r-L)}, & r \in I_2.
	\end{cases}
\end{equation}
The complete, variationally permissible global basis set $W_N = \text{span}\{\varphi_n\}_{n=1}^N$ is constructed by concatenating Equations~\eqref{eq:basis_core}, \eqref{eq:basis_tail} and \eqref{eq:hat}:
\begin{equation}\label{eq: basis}
	\varphi_n(r) =
	\begin{cases}
		\phi_{n}^1(r), & 1\leq n\leq N_1,\\[3pt]
		\phi^2_{n-N_1}(r), & N_1+1\leq n\leq N_1+N_2,\\[3pt]
		\phi^3(r),& n=N,
	\end{cases}
\end{equation}
yielding a total dimension of $N = N_1 + N_2 + 1$. 

\subsection{Variational Formulation and Numerical Implementation}
To numerically solve the RDEQ while strictly circumventing variational collapse \cite{hill1994solution}, we adopt IDOM. Rather than addressing the original Dirac operator $\mathbf{D}_r$ directly, which is plagued by a spectrum that is not bounded from below or above, we seek the stationary points of the Rayleigh quotient associated with the inverse operator $\mathbf{G} = \mathbf{D}_r^{-1}$.

To avoid the explicit construction of $\mathbf{D}_r^{-1}$, we employ the trial function transformation $\phi = \mathbf{D}_r \psi$. This reformulates the variational problem as finding the eigenpairs $(E, \psi)$ that satisfy: 
\begin{equation}\label{eq:Rayleigh}
	\frac{(\psi, \mathbf{D}_r \psi)}{(\psi, \mathbf{D}_r^2 \psi)} = \frac{1}{E}.
\end{equation}
We define the discrete approximation space for the two-component radial spinor as $\boldsymbol{V}_N = W_N \times W_N$, where $W_N$ is the global basis constructed in Equation~\eqref{eq: basis}. The Galerkin weak form of Equation~\eqref{eq:Rayleigh} leads to the following generalized eigenvalue problem: find $(E_N, \Psi_N) \in \mathbb{R} \times \boldsymbol{V}_N$ such that for all test functions $\Phi \in \boldsymbol{V}_N$,
\begin{equation}\label{eq:IDOM_Galerkin}
	\mathcal{B}(\Psi_N, \Phi) = \frac{1}{E_N} \mathcal{A}(\Psi_N, \Phi),
\end{equation}
where the bilinear forms are defined utilizing the standard Lebesgue measure on $(0, \infty)$:
\begin{align*}
	\mathcal{B}(\Psi, \Phi) = \int_0^\infty \Phi^\dagger \mathbf{D}_r \Psi ~ \mathrm dr, \quad 
	\mathcal{A}(\Psi, \Phi) = \int_0^\infty (\mathbf{D}_r\Phi)^\dagger (\mathbf{D}_r \Psi) ~ \mathrm dr.
\end{align*}
Substituting the basis expansion yields the matrix equation $\mathbf{B}\mathbf{c} = \frac{1}{E} \mathbf{A}\mathbf{c}$. 

\subsubsection*{Conditioning and Algebraic Solvers}
It is crucial to address the numerical stability and conditioning of this formulation. Because the Galerkin weak form utilizes the standard Lebesgue measure $\mathrm{d}r$, whereas GLOFs are strictly orthogonal only under a weighted measure, the global basis is non-orthogonal. Consequently, the resulting Gram and stiffness matrices are non-diagonal. Furthermore, the squared Dirac operator in IDOM amplifies the condition number $\kappa(\mathbf{A})$ as the basis size $N$ increases. However, the ALLSM framework maintains a highly compact basis representation: for the degrees of freedom required to reach a relative accuracy of $10^{-10}$ ($N \approx 60-120$), the condition number remains safely within the operational threshold of double-precision arithmetic.

Computationally, since the basis size is relatively small, the resulting generalized eigenvalue problem can typically be solved within one second on standard hardware using robust solvers such as the QZ algorithm. Should further performance optimization be required, the symmetric positive definite nature of matrix $\mathbf{A}$ allows for the application of Cholesky decomposition to reduce the generalized system to a standard eigenvalue problem, providing an accelerated route for spectral resolution without the need for iterative preconditioning.

\section{Results and Discussion}\label{sec: 5}
In this section, we comprehensively evaluate the accuracy, robustness, and physical applicability of the proposed adaptive spectral-element framework across a diverse hierarchy of relativistic quantum models. Unless otherwise specified, all numerical evaluations and error analyses are conducted in Hartree atomic units. 

We begin at the extreme limits of relativistic core interactions, resolving the severe non-polynomial spatial singularities inherent to heavy hydrogen-like ions (Subsection \ref{subsec: 5.1}). Transitioning to more physically realistic environments, we benchmark the framework against finite nuclear charge distributions (Subsection \ref{subsec: 5.2}), addressing a critical requirement for high-fidelity all-electron simulations of heavy elements. The applicability of our method is then extended to model plasma and condensed matter effects via the screened Yukawa potential (Subsection \ref{subsec: 5.3}). Furthermore, to rigorously demonstrate the structural robustness and generalization capability of our basis beyond singular configurations, the framework is applied to the pure confining relativistic harmonic oscillator (Subsection \ref{subsec: 5.4}). Here, we establish a direct, side-by-side numerical validation against the modern high-precision finite-element code \texttt{dftatom}. Finally, to further illustrate the versatility of our approach in handling complex effective core interactions, we evaluate the Hellmann potential (Subsection \ref{subsec: 5.5}).

\subsection{Accuracy and Convergence for Coulomb Potential}\label{subsec: 5.1}
To rigorously benchmark the numerical accuracy of the proposed framework, we first consider the RDEQ under the Coulomb potential, $V(r) = -Z/r$. This system serves as an indispensable cornerstone in relativistic quantum mechanics because its bound-state eigenvalues are analytically exact, governed by the fine-structure formula,
\begin{equation}\label{for:E}
	E_{n\kappa} = \frac{c^2}{\sqrt{1+\frac{(Z/c)^2}{(n-|\kappa|+s)^2}}} - c^2, \quad \text{where } s = \sqrt{\kappa^2-(Z/c)^2}.
\end{equation}
Here, $n \geq 1$ is the principal quantum number and $\kappa$ is the relativistic angular momentum quantum number.

As established in Section \ref{sec: 2}, the primary numerical hurdle for this potential is the severe non-polynomial singularity $r^s$ at the origin. To demonstrate the necessity and efficacy of our dual-adaptive spatial strategy, we systematically compare three distinct discretization schemes: (1) a standard Laguerre-only method (SLM) lacking singularity treatment; (2) a non-adaptive Log-Laguerre spectral method (LLSM) utilizing the composite GLOF-Laguerre basis but with a fixed scaling parameter $\beta$; and (3) our fully proposed adaptive Log-Laguerre spectral method (ALLSM).

\begin{figure}[h]
	\begin{minipage}{0.475\linewidth}
		\begin{center}
			\includegraphics[scale=0.325]{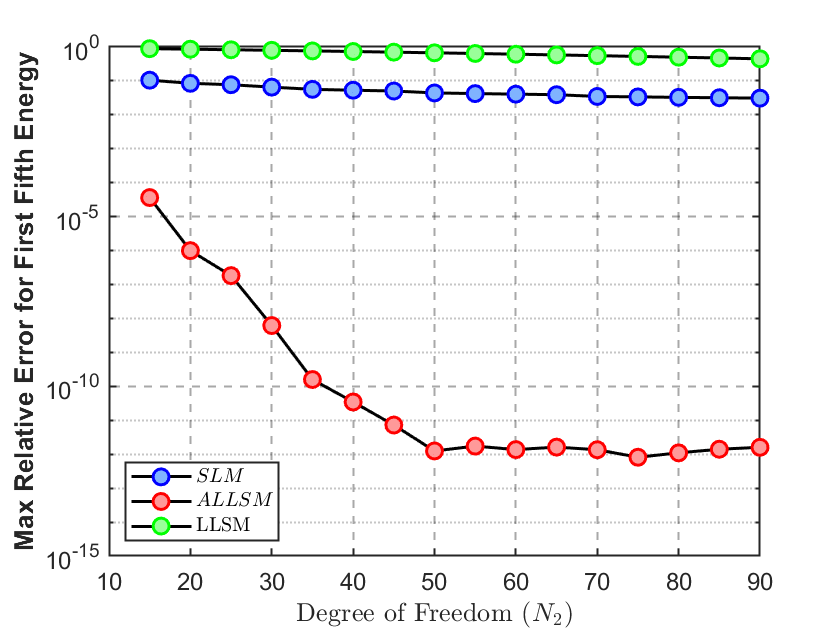}
		\end{center}
	\end{minipage}
	\begin{minipage}{0.475\linewidth}
		\begin{center}
			\includegraphics[scale=0.325]{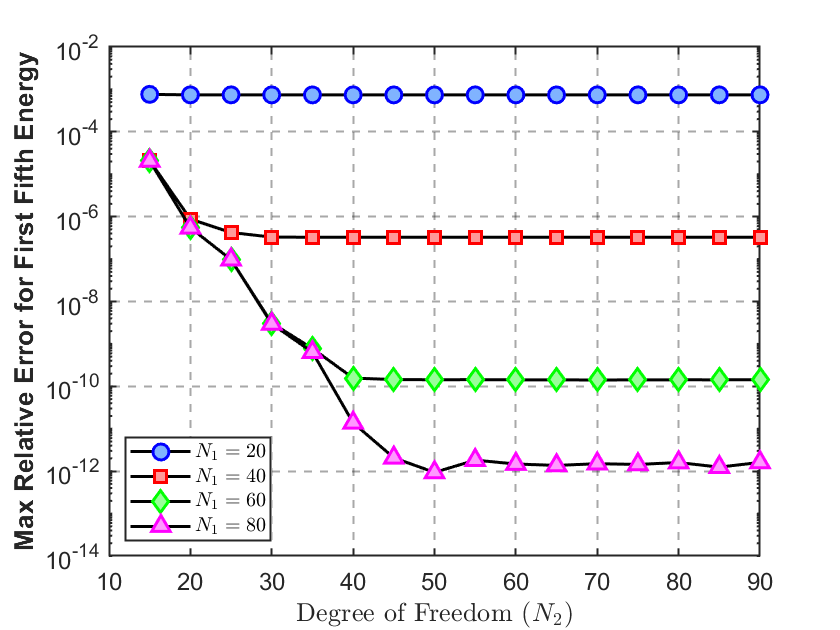}
		\end{center}
	\end{minipage}
	\caption{Example: Coulomb Potential for Uranium ion ($Z=92, \kappa=-1$). \quad \textbf{Left}: Comparison of maximum relative error versus $N_2$ for SLM, LLSM, and ALLSM, highlighting the failure of standard schemes in the highly relativistic regime. \quad\textbf{Right}: Error of ALLSM as a function of $N_2$ for various fixed $N_1$, demonstrating that our method restores strict exponential convergence by dynamically matching the state-dependent exponential tails.}
	\label{fig: 5} 
\end{figure}

We deliberately select the Uranium ion, $Z=92, \kappa=-1$ for our primary benchmark. In such a highly relativistic heavy-element regime, the singularity exponent $s$ deviates significantly from an integer, creating a pathological near-core environment that severely degrades. In left panel of Figure \ref{fig: 5}, we fix the singular GLOFs basis size at $N_1 = 80$ and plot the maximum relative error among the first five eigenvalues as a function of the Laguerre basis size, $N_2$. The results are unequivocal, the SLM fails to achieve meaningful precision, and the LLSM stagnates prematurely due to far-field basis mismatch. Instead, the proposed ALLSM dynamically perfectly matches the diverse state-dependent exponential tails, delivering superior accuracy that surpasses the alternative schemes by several orders of magnitude. The Right Panel of Figure \ref{fig: 5} further corroborates this by plotting the error against $N_2$ for various fixed $N_1$ values, clearly confirming that ALLSM successfully restores strict exponential convergence.

\begin{figure}[h]
	\begin{center}
		\includegraphics[scale=0.35]{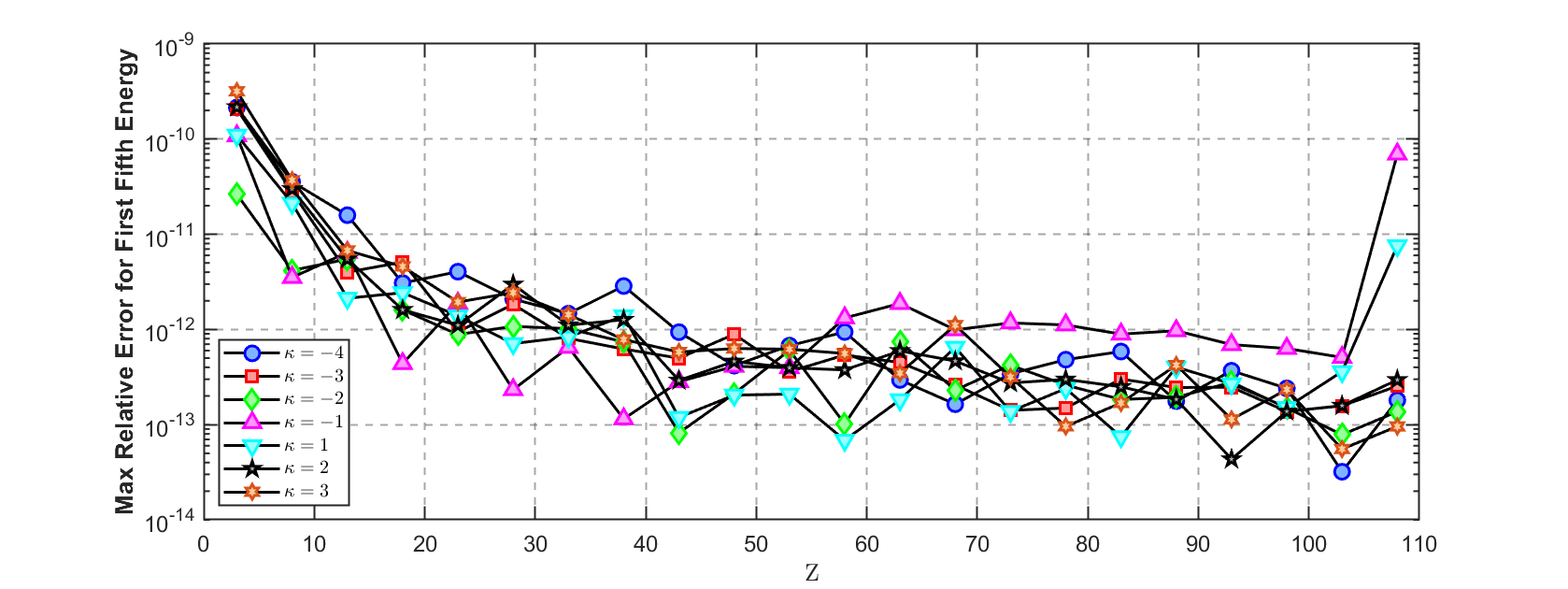}
	\end{center}
	\caption{Example: Coulomb Potential.\quad Maximum relative error of ALLSM for different $Z$ and $\kappa$ values.}
	\label{fig: 6} 
\end{figure}

To ensure these exceptional results are not an artifact of a specific parameter set, we extend the error analysis across various physical states characterized by different $Z$ and $\kappa$ combinations. As depicted in Figure \ref{fig: 6}, ALLSM consistently maintains high-fidelity precision—typically $10^{-10}$ or better—across diverse relativistic regimes. This stringent benchmarking definitively confirms the robustness of the adaptive spectral-element approach in resolving the extreme multi-scale challenges inherent to heavy-element core singularities.


\subsection{Gaussian Nuclear Model}\label{subsec: 5.2}
For highly relativistic heavy-element calculations, the idealized point-charge nuclear model becomes physically inadequate. It unphysically exaggerates the relativistic contraction of inner-shell electrons and introduces severe artifacts in core-property predictions. To provide a rigorous computational kernel for realistic all-electron simulations, we evaluate our framework using the finite nuclear charge model characterized by a Gaussian density distribution, 
\begin{equation*}
	\rho(r) = \rho_0 e^{-(r/R)^2}.
\end{equation*}
A critical physical consequence of this smeared distribution is the mathematical regularization of the strict Coulomb singularity. The resulting electrostatic potential is given by the error function, 
\begin{equation*}
	V(r) = -\frac{Z}{r} \text{erf}\left(\frac{r}{R}\right).
\end{equation*}
While this finite-nucleus potential strictly eliminates the non-polynomial $r^s$ singularity at the origin, it introduces a steep, multi-scale transition near the microscopic nuclear boundary $R$ that tightly couples with the diffuse far-field decay, presenting a severe challenge for standard fixed-scale basis sets.

\begin{figure}[h]
	\centering
	\begin{minipage}{0.475\linewidth}
		\begin{center}
			\includegraphics[scale=0.325]{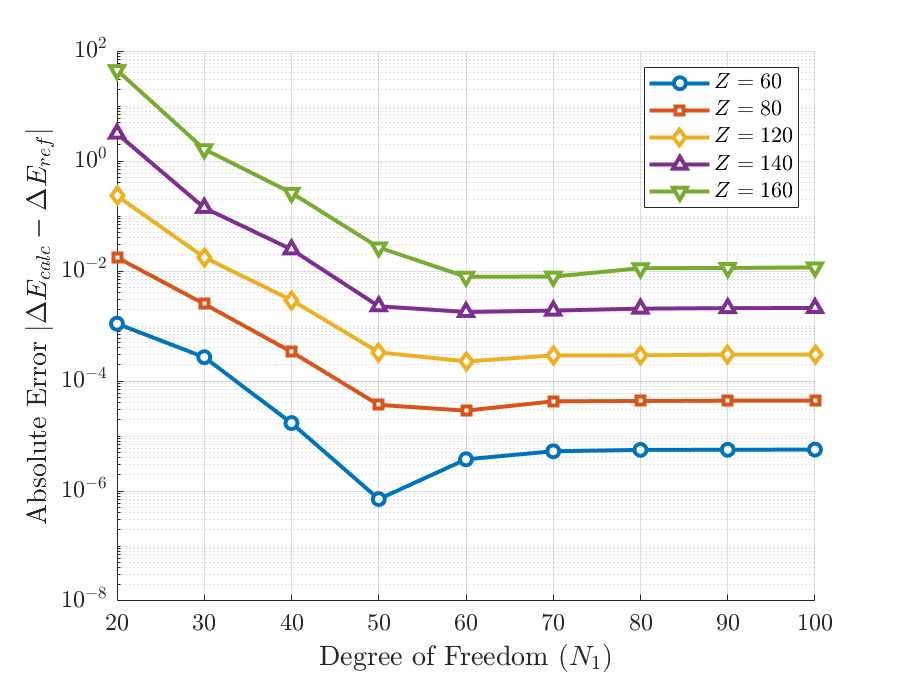}
		\end{center}
	\end{minipage}
	\begin{minipage}{0.475\linewidth}
		\begin{center}
			\includegraphics[scale=0.325]{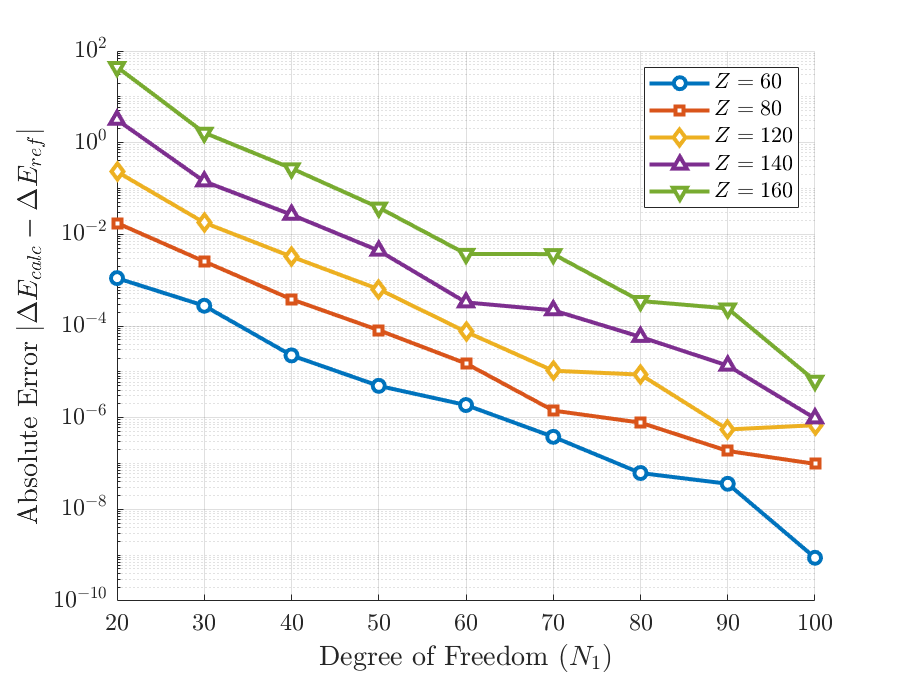}
		\end{center}
	\end{minipage}
	\caption{Example: Gaussian Nuclear Model.\quad \textbf{Left}: Comparison with established benchmarks from Andrae \textit{et al.} \cite{andrae2000finite}; the slight error increase at high basis sizes occurs as the ALLSM precision exceeds the literature's limits. \quad \textbf{Right}: Internal consistency test against an ultra-high-fidelity reference ($N_1=120, N_2=60$), confirming strict convergence.}
	\label{fig: 11}
\end{figure}

To benchmark the efficacy of the proposed ALLSM, we first compare our computed ground-state energies with the established high-precision values provided by Andrae \textit{et al.} \cite{andrae2000finite}. As depicted in the Left Panel of Figure \ref{fig: 11}, the absolute error initially plummets exponentially as the basis size increases. However, a slight artificial increase in error is observed at very high degrees of freedom. Rather than indicating an algorithmic instability, this deviation occurs precisely because our method surpasses the precision limits of the established literature benchmarks. To rigorously corroborate this intrinsic stability, we compute an ultra-high-fidelity internal reference utilizing an expanded basis with $N_1 = 120, N_2 = 60$. Against this stringent internal standard, in Right Panel of Figure \ref{fig: 11}, the error strictly and monotonically decreases down to machine precision, confirming the absolute numerical robustness and high-order convergence of our scheme.

Finally, to furnish the computational chemistry community with ultra-precise reference data, we compute the finite-nucleus energy shifts, defined as $\Delta E_{gs} = |E_{gs}^{\text{Gaussian}} - E_{gs}^{\text{Coulomb}}|$, where the point-nucleus reference energy $E_{gs}^{\text{Coulomb}}$ is exactly evaluated using the analytical fine-structure formula Equation \eqref{for:E}. These shifts quantify the massive energetic impact of nuclear smearing, which becomes exceptionally critical for superheavy elements where the point-nucleus model approaches the unphysical $Z \approx 137$ collapse catastrophe. As summarized in Table \ref{tab:results}, we successfully track these shifts across a range of heavy to superheavy elements up to the extreme relativistic limit with $Z=137$. These results not only show excellent agreement with existing models but seamlessly extend the precision boundary, serving as new high-fidelity benchmarks for future relativistic algorithm developments.

\begin{table}[ht]
	\centering
	\caption{Example: Gaussian Nuclear Model.\quad Ground-state energy shifts $\Delta E_{gs}$ for the Gaussian nuclear model. The model elegantly regularizes the point-nucleus catastrophe at $Z=137$.}
	\label{tab:results}
	\begin{tabular}{lccccc}
		\toprule 
		$Z$ & 60 & 80 & 100 & 120 & 137 \\
		\midrule
		$\Delta E_{gs}$ & 0.2299228 &
		1.9978925&
		16.4051094&
		167.3127501&
		3842.7299905 \\ 
		\bottomrule
	\end{tabular}
\end{table}

\begin{figure}[h]
	\centering
	\begin{minipage}{0.475\linewidth}
		\begin{center}
			\includegraphics[scale=0.325]{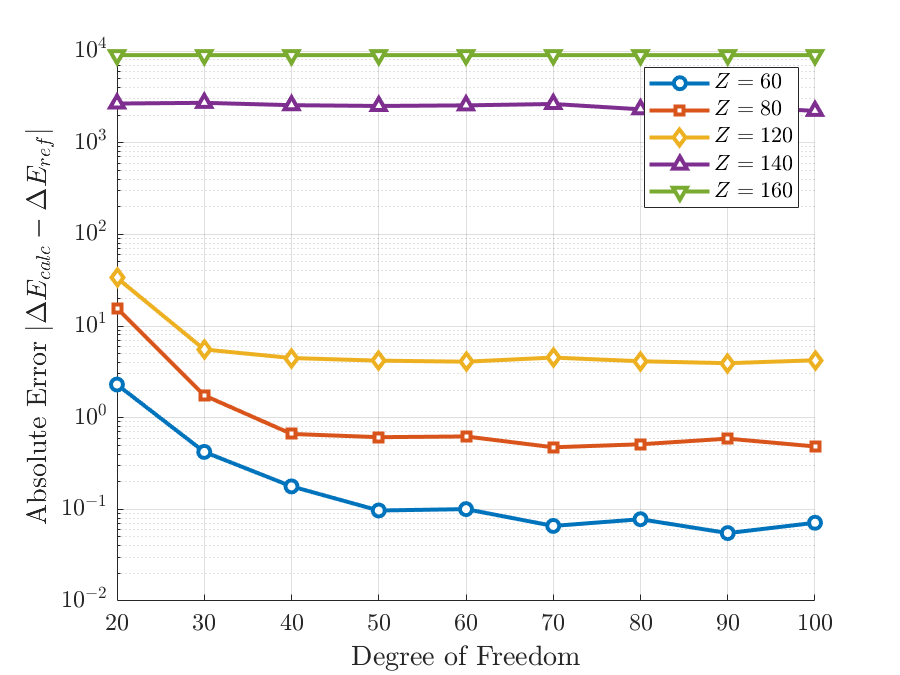}
		\end{center}
	\end{minipage}
	\begin{minipage}{0.475\linewidth}
		\begin{center}
			\includegraphics[scale=0.325]{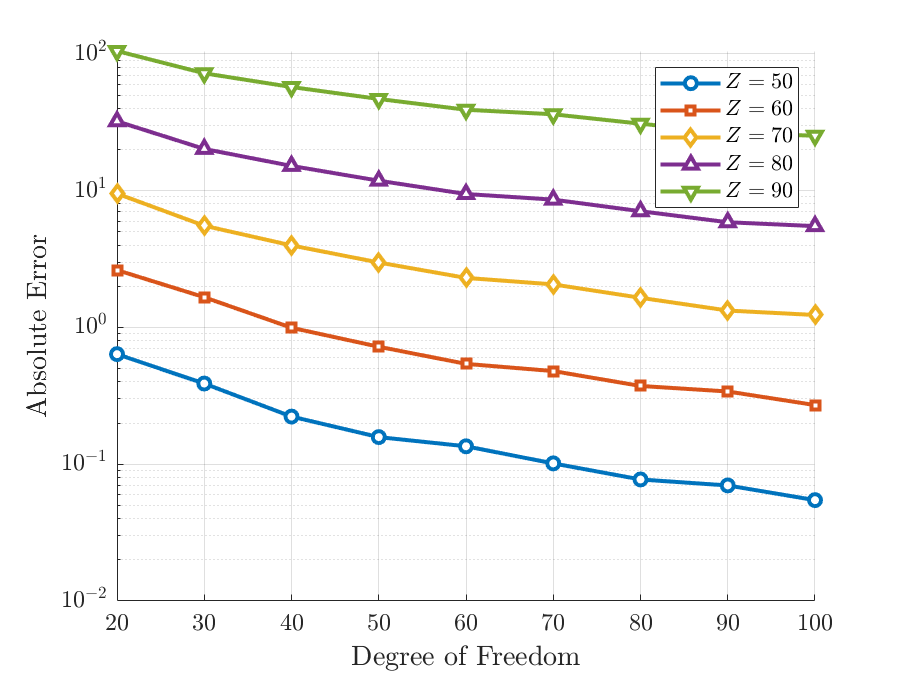}
		\end{center}
	\end{minipage}
	\caption{Convergence behavior of the numerical error for different nuclear potentials using ALSM. \quad \textbf{{Left}}: Error evolution for the Gaussian nuclear model. \quad \textbf{Right}: Error evolution for the singular Coulomb potential. A highly identical failure pattern is observed, indicating that the multi-scale steep gradients and spinor kinetic mismatch in the Gaussian model induce a "pseudo-singularity" that identically restrains the convergence of standard ALSM.}
	\label{fig: 111}
\end{figure}

A profound numerical phenomenon is observed regarding the structural limits of global polynomial bases in this model. The error-function-based potential is strictly analytic at the origin:
\begin{equation*}
	V(r) = -\frac{Z}{r} \text{erf}\left(\frac{r}{R}\right) = -\frac{2Z}{\sqrt{\pi} R} + \frac{2Z}{3\sqrt{\pi} R^3} r^2 - \frac{Z}{5\sqrt{\pi} R^5} r^4 + \mathcal{O}(r^6),
\end{equation*}
yielding a finite limiting value of $V(0) = -2Z/(\sqrt{\pi}R)$. Since the fractional $r^s$ behavior is mathematically eliminated, one might intuitively expect a standard global polynomial basis (e.g., ALSM) to achieve rapid exponential convergence. However, as illustrated in the Left Panel of Figure \ref{fig: 111}, the numerical error of ALSM exhibits a persistent, catastrophic stagnation. Strikingly, the error evolution for this non-singular Gaussian model (Left Panel) identically mirrors the failure pattern observed in the singular Coulomb potential (Right Panel).

This unexpected ``pseudo-singularity'' is rooted in the kinetic coupling of the Dirac operator. For heavy elements, the microscopic nuclear radius, $R \sim 10^{-5}$ a.u. creates an extremely localized, steep gradient boundary layer. Furthermore, approximating both the upper and lower spinor components within the identical global spatial basis induces a severe kinetic balance mismatch near this abrupt structural transition. When utilizing operator reformulations, such as IDOM, this rigidity prevents the simultaneous resolution of both spinor dynamics, severely restraining the convergence of pure Laguerre functions. Consequently, the computational barrier is not the non-analytic power law, but rather the multi-scale gradient and the associated dynamical mismatch. 

In this regime, ALLSM provides an elegant topological remedy. The logarithmic mapping inherent to GLOFs acts as an intrinsic infinite mesh refinement mechanism near the origin. This smoothly absorbs the steep boundary layer and provides the spatial flexibility necessary to resolve the upper-lower spinor mismatch. Thus, even in systems completely devoid of traditional singularities, ALLSM maintains remarkable robustness, significantly outperforming standard spectral frameworks by neutralizing these implicit structural pathologies.

\subsection{Yukawa Potential}\label{subsec: 5.3}
To demonstrate the applicability of our method to physically relevant environments, such as atoms in plasmas or electrons in condensed matter, we consider the screened Yukawa potential:
\begin{equation*}
	V(r) = -V_0 \frac{e^{-\lambda r}}{r},
\end{equation*}
where $\lambda$ is the screening parameter and $V_0$ is the coupling strength. As $\lambda$ increases, the effective interaction range decreases, eventually leading to the disappearance of bound states near the continuum threshold. This transition region is numerically sensitive and prone to variational instabilities.

Using ALLSM, we compute the ground-state energy for two benchmark cases adopted from previous literature: Case 1 with $V_0 = 0.1c, \lambda = 0.01c$ \cite{krauthauser2002basis}, and Case 2 with $V_0 = 0.4c, \lambda = 0.07c$ \cite{baye2014accurate}. To rigorously validate our framework against established numerical standards, Table \ref{tab:yukawa_benchmark} presents a direct, side-by-side numerical comparison of our computed energies against the highly accurate reference data obtained via B-spline \cite{krauthauser2002basis} and Lagrange mesh methods \cite{baye2014accurate}. The ALLSM results strictly match the literature benchmarks up to 10-12 digits, explicitly confirming the extreme precision of our spectral-element implementation.

\begin{table}[ht]
	\centering
	\caption{Direct numerical comparison of the ground-state energy for the Yukawa potential. The ALLSM results are benchmarked against established high-precision data from B-spline and Lagrange mesh methodologies.}
	\label{tab:yukawa_benchmark}
	\small
	\begin{tabular}{l l l l}
		\toprule
		{Parameter Set} & {Method} & {$E_{gs}$ (Hartree)} & {Reference} \\
		\midrule
		\multirow{2}{*}{Case 1: $V_0 = 0.1c, \lambda = 0.01c$} 
		& ALLSM & \textbf{-76.672555928755170} & This work \\
		& Lagrange Mesh & -76.672555928706560 & Ref. \cite{krauthauser2002basis} \\
		\midrule
		\multirow{2}{*}{Case 2: $V_0 = 0.4c, \lambda = 0.07c$} 
		& ALLSM & \textbf{-4.859949240857070$\times 10^{3}$} & This work \\
		& Lagrange Mesh & -4.859949240864677e+03 & Ref. \cite{baye2014accurate} \\
		\bottomrule
	\end{tabular}
\end{table}

\begin{figure}[h]
	\begin{minipage}{0.475\linewidth}
		\begin{center}
			\includegraphics[scale=0.325]{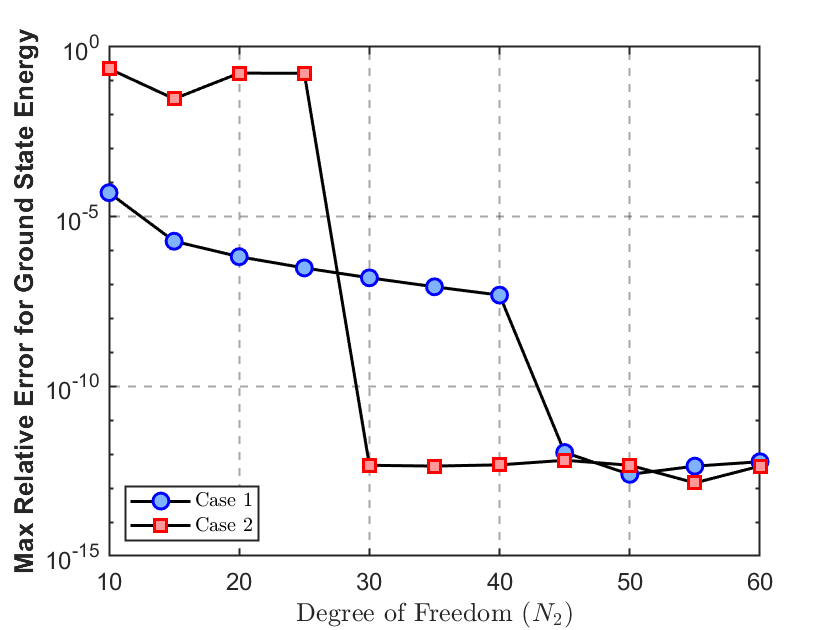}
		\end{center}
	\end{minipage}
	\begin{minipage}{0.475\linewidth}
		\begin{center}
			\includegraphics[scale=0.325]{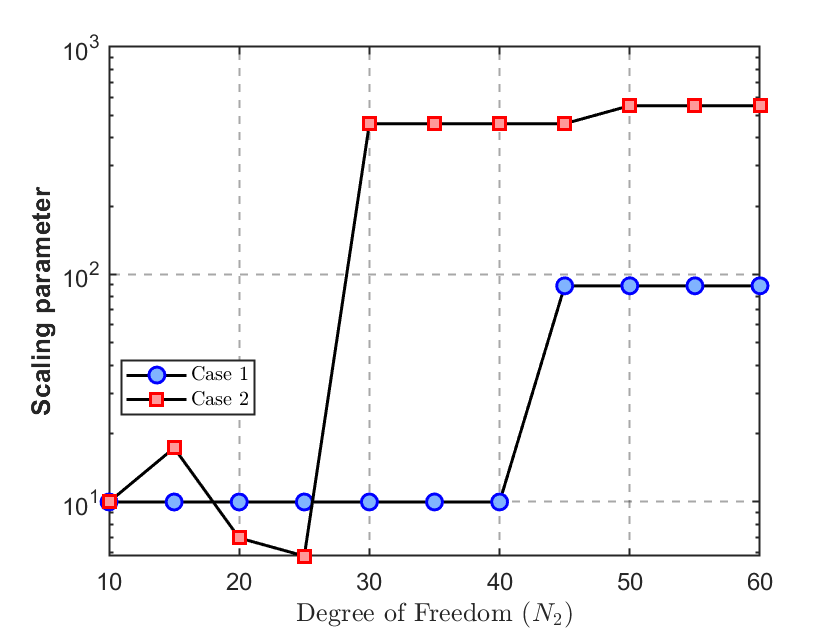}
		\end{center}
	\end{minipage}
	\caption{Example: Yukawa Potential.\quad  \textbf{Left}: Maximum relative error of the ground-state energy for Case 1 ($V_0 = 0.1c, \lambda = 0.01c$) and Case 2 ($V_0 = 0.4c, \lambda = 0.07c$). The results demonstrate that ALLSM maintains high precision across different screening regimes. \quad \textbf{Right}: Evolution of the adaptive scaling parameter as a function of the degrees of freedom $N_2$. The scaling factor effectively converges from a sub-optimal starting point to its appropriate value as the basis size increases, ensuring a minimized error.}
	\label{fig: 7}
\end{figure}

Beyond isolated energy validation, the internal convergence mechanism of ALLSM is displayed in the Left Panel of Figure \ref{fig: 7}. The method successfully maintains exponential convergence for these screened potentials. The Right Panel of Figure \ref{fig: 7} illustrates the behavior of the adaptive scaling parameter $\beta$. At lower degrees of freedom, the optimization landscape is relatively flat; however, as the basis size $N$ increases, the algorithm robustly locks onto the global optimum, ensuring structural fidelity.

\begin{figure}[h]
	\begin{minipage}{0.475\linewidth}
		\begin{center}
			\includegraphics[scale=0.325]{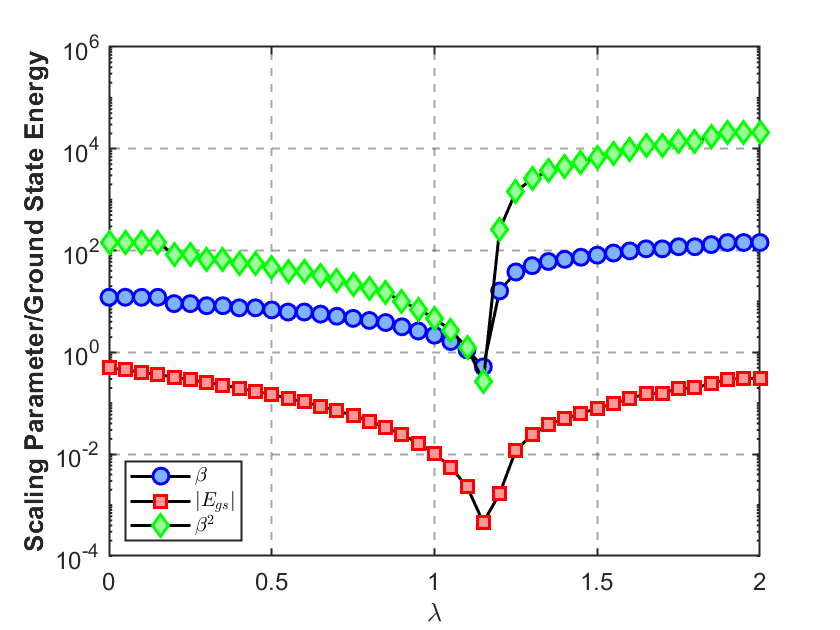}
		\end{center}
	\end{minipage}
	\begin{minipage}{0.475\linewidth}
		\begin{center}
			\includegraphics[scale=0.325]{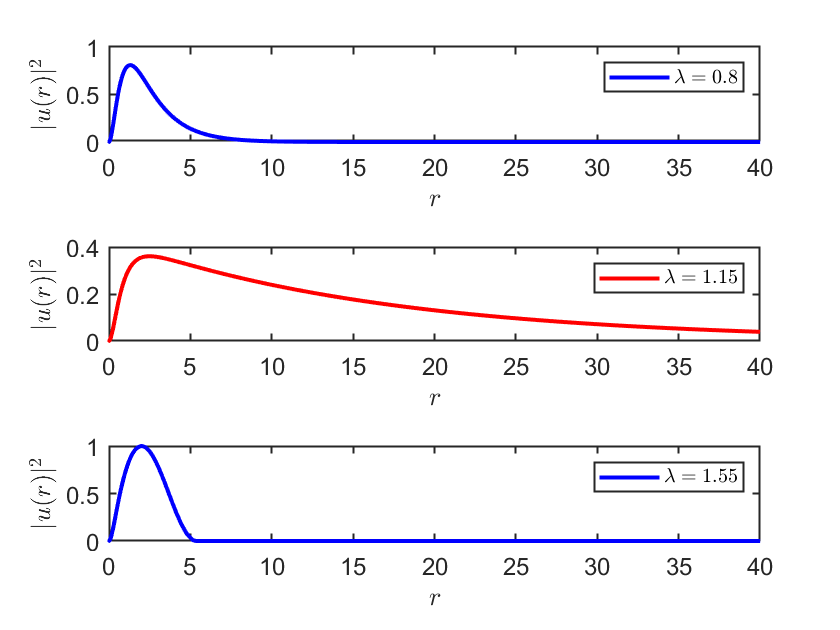}
		\end{center}
	\end{minipage}
	\caption{Example: Yukawa Potential. \textbf{Left}: Comparative evolution of the ground-state energy $E_{\text{gs}}$ and the optimal scaling parameters, $\beta$ and $\beta^2$, as a function of the screening parameter $\lambda$. The synchronized decaying trend suggests that the adaptive scaling factor effectively tracks the energy shift as the system approaches the continuum. 
	\textbf{Right}: Radial wavefunctions for various $\lambda$ values approaching $\lambda_{\text{crit}}$. As the binding energy diminishes, the wavefunctions exhibit a transition from localized states to increasingly diffuse and uniform distributions, reflecting the loss of confinement near the ionization threshold.}
	\label{fig: 8}
\end{figure}

Finally, we investigate the critical screening parameter, $\lambda_{\text{crit}}$, at which the ground state energy vanishes into the continuum. Setting $V_0 = 1$, we identify the critical value at approximately $\lambda_{\text{crit}} \approx 1.190613$. In the Left Panel of Figure \ref{fig: 8}, we plot the ground state energy $E_{gs}$ and the optimal scaling parameter $\beta$ (specifically $\beta^2$) as functions of $\lambda$. It is mathematically evident that $E_{gs}$ and $\beta^2$ exhibit a highly synchronized decaying trend. 
This strong correlation is rooted in the physical asymptotics of the wavefunction: the scaling parameter $\beta$ actively tracks the characteristic decay length of the system, which is intrinsically governed by $\sqrt{c^2 - (E/c)^2}$ for the dominant far-field tail. This implies that the optimized $\beta$ can serve as a reliable intrinsic proxy for monitoring the effective spatial extent of the wavefunction near the continuum threshold. The Right Panel of Figure \ref{fig: 8} corroborates this physics: as $\lambda$ increases towards the critical ionization threshold, the binding energy diminishes, and the corresponding wavefunctions become dramatically more diffuse and spatially extended, reflecting the loss of electrostatic confinement.

\subsection{Harmonic Oscillator Potential}\label{subsec: 5.4}
To further demonstrate the generalization performance and structural robustness of the ALLSM framework, we evaluate the relativistic harmonic oscillator, defined by the pure scalar confining potential $V(r) = \frac{1}{2}r^2$. 

Unlike the Coulombic systems, the harmonic oscillator potential does not exhibit a non-polynomial $r^s$ singularity at the origin, and its wavefunctions do not exhibit the standard exponential decay ($\sim e^{-\lambda r}$). This makes it an exceptionally rigorous test case for our framework: it evaluates whether the Log-Laguerre and adaptive Laguerre basis sets can flexibly seamlessly adapt to non-singular, purely confined systems without suffering from numerical over-completeness or performance degradation.

Because simple analytical solutions for this relativistic model do not exist, numerical benchmarking against established, independent codes is essential. To this end, we directly compare our ALLSM results against \texttt{dftatom} \cite{vcertik2013dftatom}, a modern, high-precision finite-element atomic code specifically highlighted for its rigorous boundary convergence and robust mesh grading capabilities. 

As detailed in Table \ref{table: kap}, the eigenvalues computed by ALLSM exhibit outstanding agreement with the \texttt{dftatom} reference data across a wide range of quantum states. The ALLSM framework handles this purely confining potential seamlessly, achieving extreme precision using an extremely compact basis set. This side-by-side validation explicitly confirms that our proposed spectral-element methodology is not strictly confined to singularly-peaked potentials, but serves as a universally highly accurate, pollution-free numerical engine for diverse relativistic physical models.

\begin{table}[!ht]
	\centering
	\caption{Example: Harmonic Oscillator Potential. The eigenvalues computed with different principal quantum numbers $n$ and relativistic angular momentum quantum numbers $\kappa$. The ALLSM results are benchmarked against high-precision reference data from \texttt{dftatom} \cite{vcertik2013dftatom}.}
	\label{table: kap}
	\small
	\begin{tabular}{c c l l}
		\toprule
		$n$ & $\kappa$ & \texttt{dftatom} Reference & ALLSM Energy \\ 
		\midrule
		$1$ & $-1$ & 1.49999501 & 1.499995007841790 \\ 
		\midrule
		\multirow{3}{*}{$2$} 
		& $-1$ & 3.49989517 & 3.499895170505624 \\ 
		& $-2$ & 2.49997504 & 2.499975038957928 \\ 
		& $1$  & 2.49993510 & 2.499935104453471 \\ 
		\midrule
		\multirow{5}{*}{$3$} 
		& $-1$ & 5.49971548 & 5.499715477584687 \\ 
		& $-2$ & 4.49983527 & 4.499835270624317 \\ 
		& $1$  & 4.49979534 & 4.499795339215780 \\ 
		& $-3$ & 3.49994176 & 3.499941758258501 \\ 
		& $2$  & 3.49987520 & 3.499875203280681 \\ 
		\midrule
		\multirow{7}{*}{$4$} 
		& $-1$ & 7.49945594 & 7.499455941229826 \\ 
		& $-2$ & 6.49961565 & 6.499615652712237 \\ 
		& $1$  & 6.49957572 & 6.499575724505121 \\ 
		& $-3$ & 5.49976206 & 5.499762061157526 \\ 
		& $2$  & 5.49969551 & 5.499695511567552 \\ 
		& $-4$ & 4.49989517 & 4.499895166587521 \\ 
		& $3$  & 4.49980199 & 4.499801993515575 \\ 
		\midrule
		\multirow{7}{*}{$5$} 
		& $-1$ & 9.49911657 & 9.499116573671927 \\ 
		& $-2$ & 8.49931620 & 8.499316197503504 \\ 
		& $1$  & 8.49927627 & 8.499276272519637 \\ 
		& $-3$ & 7.49950252 & 7.499502520728129 \\ 
		& $2$  & 7.49943598 & 7.499435976427776 \\ 
		& $-4$ & 6.49967554 & 6.499675647402910 \\ 
		& $3$  & 6.49958238 & 6.499582482058031 \\ 
		\midrule
		\multirow{7}{*}{$6$} 
		& $-1$ & 11.49869739 & 11.498697387101856 \\ 
		& $-2$ & 10.49893692 & 10.498936917203537 \\ 
		& $1$  & 10.49889699 & 10.498896995446557 \\ 
		& $-3$ & 9.49916315 & 9.499163149106607 \\ 
		& $2$  & 9.49909661 & 9.499096610030392 \\ 
		& $-4$ & 8.49937608 & 8.499380441662652 \\ 
		& $3$  & 8.49928292 & 8.499287401344190 \\ 
		\midrule
		\multirow{7}{*}{$7$} 
		& $-1$ & 13.49819839 & 13.498198393761413 \\ 
		& $-2$ & 12.49847782 & 12.498477824050497 \\ 
		& $1$  & 12.49843790 & 12.498437905396713 \\ 
		& $-3$ & 11.49874396 & 11.498743958352861 \\ 
		& $2$  & 11.49867742 & 11.498677424722700 \\ 
		& $-4$ & 10.49899680 & 10.499046796689072 \\ 
		& $3$  & 10.49890364 & 10.498954300546757 \\ 
		\bottomrule
	\end{tabular}
\end{table}

\subsection{Hellmann Potential}\label{subsec: 5.5}
To further demonstrate the versatility of our adaptive framework in handling complex, multi-term effective interactions, we investigate the Hellmann potential. Widely utilized in atomic physics and quantum chemistry to model electron-ion interactions and construct empirical pseudopotentials \cite{korsch1982milne,hamzavi2012approximate}, the Hellmann potential is formulated as the superposition of a long-range Coulomb attraction and a short-range screened Yukawa interaction:
\begin{equation*}
	V(r) = -\frac{Z}{r} + V_0\frac{e^{-\lambda r}}{r},
\end{equation*}
where $Z$ dictates the bare Coulomb strength, $V_0$ represents the coupling amplitude of the Yukawa core, and $\lambda$ acts as the environmental screening parameter.

The physical landscape of this potential is intricately governed by the delicate balance between these two terms, exhibiting significant sensitivity to $\lambda$. When $V_0 > 0$, the Yukawa component introduces a short-range repulsion. In this regime, an increase in $\lambda$ rapidly suppresses the repulsive tail, leading to a deeper effective potential well and consequently more strongly bound states. Conversely, when $V_0 < 0$, the Yukawa term provides a supplementary short-range attraction. Here, increasing the screening $\lambda$ curtails this additional attractive range, thereby shallowing the well.

Despite its critical utility in chemical physics, rigorous, high-precision numerical benchmarks for the Hellmann potential within the unsimplified RDEQ framework remain notably absent from the current literature. To bridge this critical data gap, we deploy ALLSM to compute extensive reference eigenvalues across a diverse matrix of parameter configurations, as cataloged in Table \ref{tab:hellmann_benchmark1}.

ALLSM perfectly resolves the subtle energetic shifts induced by these competing interactions. For instance, considering the supplementary attraction regime with $Z=2, V_0=-1$, our data in Table \ref{tab:hellmann_benchmark1} clearly illustrate that as the screening parameter $\lambda$ incrementally increases from $0.001$ to $0.01$, the absolute binding energy strictly decreases. This ultra-precise numerical tracking perfectly corroborates the physical intuition, namely, a heightened screening parameter dynamically truncates the effective range of the Yukawa attraction, systematically ``lifting'' the tightly bound energy levels toward the continuum threshold.
\begin{table}[ht]
	\centering
	\caption{Example: Hellmann potential.\quad The different state energy for Hellmann potential.}
	\label{tab:hellmann_benchmark1}
	\begin{tabular}{lcccc}
		\hline
		State & $\lambda$ & $Z=3,~V_0=0$ & $Z=0,~V_0=-3$ & $Z=2,~V_0=-1$ \\
		\hline
		$n=1 ,\kappa=-1$ & 0.001 & -4.5005393 & -4.4975400 & -4.4995395 \\
		& 0.005 & -4.5005393 & -4.4855580 & -4.4955455 \\
		& 0.01  & -4.5005393 & -4.4706141 & -4.4905642 \\
		\hline
		$n=2 ,\kappa=-1$ & 0.001 & -1.1251685 & -1.1221715 & -1.1241695 \\
		& 0.005 & -1.1251685 & -1.1102432 & -1.1201934 \\
		& 0.01  & -1.1251685 & -1.0954662 & -1.1152677 \\
		\hline
		$n=3 ,\kappa=-1$ & 0.001 & -0.5000599 & -0.4970667 & -0.4990622 \\
		& 0.005 & -0.5000599 & -0.4852272 & -0.4951157 \\
		& 0.01  & -0.5000599 & -0.4707237 & -0.4902811 \\
		\hline
		$n=2 ,\kappa=-2$ & 0.001 & -1.1250337 & -1.1220362 & -1.1240345 \\
		& 0.005 & -1.1250337 & -1.1100960 & -1.1200545 \\
		& 0.01  & -1.1250337 & -1.0952820 & -1.1151165 \\
		\hline
		$n=3 ,\kappa=-2$ & 0.001 & -0.5000200 & -0.4970262 & -0.4990221 \\
		& 0.005 & -0.5000200 & -0.4851750 & -0.4950716 \\
		& 0.01  & -0.5000200 & -0.4706352 & -0.4902250 \\
		\hline
		$n=2 ,\kappa=1$ & 0.001 & -1.1251685 & -1.1221710 & -1.1241694 \\
		& 0.005 & -1.1251685 & -1.1102308 & -1.1201893 \\
		& 0.01  & -1.1251685 & -1.0954168 & -1.1152513 \\
		\hline
		$n=3 ,\kappa=1$ & 0.001 & -0.5000599 & -0.4970662 & -0.4990620 \\
		& 0.005 & -0.5000599 & -0.4852149 & -0.4951116 \\
		& 0.01  & -0.5000599 & -0.4706751 & -0.4902650 \\
		\hline
		$n=3 ,\kappa=-3$ & 0.001 & -0.5000067 & -0.4970119 & -0.4990084 \\
		& 0.005 & -0.5000067 & -0.4851370 & -0.4950501 \\
		& 0.01  & -0.5000067 & -0.4705249 & -0.4901794 \\
		\hline
		$n=3 ,\kappa=2$ & 0.001 & -0.5000200 & -0.4970252 & -0.4990217 \\
		& 0.005 & -0.5000200 & -0.4851504 & -0.4950634 \\
		& 0.01  & -0.5000200 & -0.4705382 & -0.4901927 \\
		\hline
	\end{tabular}
\end{table}

\section{Conclusion}\label{sec: 6}
In this work, we have developed and rigorously validated a comprehensive, high-order adaptive spectral-element framework tailored for relativistic radial quantum eigenvalue problems. By seamlessly integrating IDOM to structurally guarantee spectral stability, our approach systematically overcomes the two most formidable approximation bottlenecks in relativistic electronic structure calculations: the multi-scale, state-dependent asymptotic decay of wavefunctions, and the non-polynomial core singularities (or steep near-core gradients).

To intrinsically resolve the pathological near-nuclear $r^s$ singularities characteristic of heavy-element Coulomb-like potentials, as well as the steep multi-scale boundary layers in finite-nucleus models, we deployed GLOFs. Crucially, this specialized basis operates as a robust ``black-box'' solver, exponentially resolving complex fractional power-law expansions without requiring any \textit{a priori} analytical extraction or specification of the singularity exponents. Concurrently, to capture the diffuse and highly oscillatory far-field tails, we introduced an Adaptive Laguerre Method framework. Driven by a dynamic, frequency-based indicator, this mechanism autonomously optimizes the characteristic length scale of the basis to perfectly envelop the physical support of the target state. This dual-adaptive spatial strategy intrinsically bypasses the need for artificial domain truncation or dense localized mesh refinement. It minimizes the required degrees of freedom, dramatically reducing computational overhead while preserving strict exponential convergence. It should be emphasized that our framework exhibits highly robust performance across both singular (non-polynomial) and purely confining (polynomial) wavefunction behaviors.

Beyond high-fidelity algorithmic development, this study provides fundamental theoretical insights into the fundamental mechanics of spectral pollution in Dirac discretizations. We demonstrated that artificial finite-interval truncation introduces severe efficiency bottlenecks, while upper-lower spinor degree of freedom imbalances act as primary vectors for variational collapse. Furthermore, we mathematically justified the absolute stabilizing supremacy of IDOM over conventional regularization schemes; IDOM topologically segregates the spectrum by collapsing the essential convex hull, ensuring a pristine, automated, and pollution-free eigenvalue resolution.

Extensive numerical benchmarking across a broad spectrum of physically relevant models—including bare Coulomb, finite-nucleus Gaussian, screened Yukawa, and the pure confining relativistic harmonic oscillator—confirms the strong robustness of the proposed framework. Crucially, the accuracy of our methodology was rigorously validated through direct, side-by-side numerical comparisons with established high-precision literature and modern finite-element benchmarks. Consistently achieving relative accuracies on the order of $10^{-10}$ across diverse relativistic regimes, this methodology provides a high-accuracy benchmark for atomic calculations. Ultimately, this work delivers a highly reliable, mathematically rigorous computational engine capable of generating ultra-precise reference data, laying a formidable mathematical foundation for the next generation of large-scale relativistic quantum chemistry simulations.

Looking forward, the mathematical strategies established in this work open several highly promising avenues for advanced quantum chemical methodology. First, the remarkable capability of GLOFs to act as a universal, ``black-box'' solver for unknown or fractional singularities extends far beyond the Radial Dirac Equation. This leads to a compelling dual perspective on the nature of singularities: while a rigorous analytical determination of the singularity index $s$ remains an intellectually stimulating mathematical challenge for future exploration, the present ALLSM effectively renders such prior knowledge optional. By demonstrating robust convergence even in the absence of an explicit $s$, our approach proves that the ``unknown'' is no longer a barrier to precision. Consequently, this positions GLOFs as a potent tool for other complex domains plagued by analytical singularities, such as addressing the electron-nucleus cusp conditions in explicitly correlated methods. Furthermore, the frequency-based adaptive scaling mechanism introduces a transformative paradigm for conventional molecular basis set construction. This intrinsic reduction in degrees of freedom will significantly alleviate the computational bottleneck of multi-center integral evaluations, offering a highly efficient computational strategy for tackling high-dimensional, many-body relativistic electronic structure problems.

\section*{Data and Code Availability}
To ensure the rigorous reproducibility of the computational benchmarks presented in this study, the core datasets—including the convergence profiles and the high-precision reference eigenvalues for the evaluated physical potentials (Coulomb, Gaussian, Yukawa, Harmonic Oscillator, and Hellmann)—are fully provided within the article's tables and figures. The core scripts (MATLAB) necessary to assemble the localized basis matrices and reproduce the generalized eigenvalue problems for the reported test cases are available from the corresponding author upon reasonable request.

\section{Acknowledgments}
The authors thank the National Natural Science Foundation of China (Grant Nos. 12471341, 12325112, and 12288101) for financial support. The authors thank the Laboratory of Mathematics and Applied Mathematics (LMAM) and the High-performance Computing Platform of Peking University for providing a stimulating research environment and office space to Sheng Chen and Shuai Wu during their visit.

\printbibliography

@book{zabloudil2005electron,
	title={Electron Scattering in Solid Matter: a theoretical and computational treatise},
	author={Zabloudil, Jan and Hammerling, Robert and Weinberger, Peter and Szunyogh, Laszlo},
	year={2005},
	publisher={Springer}
}

@article{lewin2010spectral,
	title={Spectral pollution and how to avoid it},
	author={Lewin, Mathieu and S{\'e}r{\'e}, {\'E}ric},
	journal={Proceedings of the London mathematical society},
	volume={100},
	number={3},
	pages={864--900},
	year={2010},
	publisher={Oxford University Press}
}

@article{chen2025efficient,
	title={Efficient Spectral-Element Methods in Polar Coordinates for Complex Geometries with Piecewise-Smooth Boundaries},
	author={Chen, Sheng and Sun, Weiwei and Wu, Shuai},
	journal={SIAM Journal on Scientific Computing},
	volume={47},
	number={3},
	pages={A1907--A1936},
	year={2025},
	publisher={SIAM}
}

@article{xia2021efficient,
	title={Efficient scaling and moving techniques for spectral methods in unbounded domains},
	author={Xia, Mingtao and Shao, Sihong and Chou, Tom},
	journal={SIAM Journal on Scientific Computing},
	volume={43},
	number={5},
	pages={A3244--A3268},
	year={2021},
	publisher={SIAM}
}

@article{kutzelnigg1984basis,
	title={Basis set expansion of the Dirac operator without variational collapse},
	author={Kutzelnigg, Werner},
	journal={International Journal of Quantum Chemistry},
	volume={25},
	number={1},
	pages={107--129},
	year={1984},
	publisher={Wiley Online Library}
}

@article{mark1982new,
	title={New Representation of the $\alpha$→{\textperiodcentered} p→ Operator in the Solution of Dirac-Type Equations by the Linear-Expansion Method},
	author={Mark, F and Schwarz, WHE},
	journal={Physical Review Letters},
	volume={48},
	number={10},
	pages={673},
	year={1982},
	publisher={APS}
}

@article{hill1994solution,
	title={A solution to the problem of variational collapse for the one-particle Dirac equation},
	author={Hill, Robert Nyden and Krauthauser, Carl},
	journal={Physical review letters},
	volume={72},
	number={14},
	pages={2151},
	year={1994},
	publisher={APS}
}

@article{baye2014accurate,
	title={Accurate solution of the Dirac equation on Lagrange meshes},
	author={Baye, Daniel and Filippin, Livio and Godefroid, Michel},
	journal={Physical Review E},
	volume={89},
	number={4},
	pages={043305},
	year={2014},
	publisher={APS}
}

@article{krauthauser2002basis,
	title={Basis-set methods for the Dirac equation},
	author={Krauthauser, Carl and Hill, Robert Nyden},
	journal={Canadian Journal of Physics},
	volume={80},
	number={3},
	pages={181--265},
	year={2002},
	publisher={NRC Research Press Ottawa, Canada}
}

@article{korsch1982milne,
	title={Milne's differential equation and numerical solutions of the Schrodinger equation. II. Complex energy resonance states},
	author={Korsch, HJ and Laurent, H and Mohlenkamp, R},
	journal={Journal of Physics B: Atomic and Molecular Physics},
	volume={15},
	number={1},
	pages={1},
	year={1982},
	publisher={IOP Publishing}
}

@article{andrae2000finite,
	title={Finite nuclear charge density distributions in electronic structure calculations for atoms and molecules},
	author={Andrae, Dirk},
	journal={Physics Reports},
	volume={336},
	number={6},
	pages={413--525},
	year={2000},
	publisher={Elsevier}
}

@article{hamzavi2012approximate,
	title={Approximate analytical solution of the Yukawa potential with arbitrary angular momenta},
	author={Hamzavi, M and Movahedi, M and Thylwe, K-E and Rajabi, AA},
	journal={Chinese Physics Letters},
	volume={29},
	number={8},
	pages={080302},
	year={2012},
	publisher={IOP Publishing}
}

@article{shen2009some,
	title={Some recent advances on spectral methods for unbounded domains},
	author={Shen, Jie and Wang, Li-Lian},
	journal={Communications in computational physics},
	volume={5},
	number={2-4},
	pages={195--241},
	year={2009}
}

@article{hou2007computing,
	title={Computing nearly singular solutions using pseudo-spectral methods},
	author={Hou, Thomas Y and Li, Ruo},
	journal={Journal of Computational Physics},
	volume={226},
	number={1},
	pages={379--397},
	year={2007},
	publisher={Elsevier}
}

@article{vcertik2024high,
	title={High-order finite element method for atomic structure calculations},
	author={{\v{C}}ert{\'\i}k, Ond{\v{r}}ej and Pask, John E and Fernando, Isuru and Goswami, Rohit and Sukumar, N and Collins, Lee A and Manzini, Gianmarco and Vack{\'a}{\v{r}}, Ji{\v{r}}{\'\i}},
	journal={Computer Physics Communications},
	volume={297},
	pages={109051},
	year={2024},
	publisher={Elsevier}
}

@article{schwarz1982mol,
	title={Mol},
	author={Schwarz, WHE and Wallmeier, H},
	journal={Phys.},
	volume={46},
	pages={1045--1061},
	year={1982}
}

@article{schwarz1982two,
	title={The two problems connected with Dirac-Breit-Roothaan calculations},
	author={Schwarz, WHE and Wechsel-Trakowski, E},
	journal={Chemical Physics Letters},
	volume={85},
	number={1},
	pages={94--97},
	year={1982},
	publisher={Elsevier}
}

@article{chen2022log,
	title={Log orthogonal functions: approximation properties and applications},
	author={Chen, Sheng and Shen, Jie},
	journal={IMA Journal of Numerical Analysis},
	volume={42},
	number={1},
	pages={712--743},
	year={2022},
	publisher={Oxford University Press}
}

@article{chen2020spectrally,
	title={A spectrally accurate approximation to subdiffusion equations using the log orthogonal functions},
	author={Chen, Sheng and Shen, Jie and Zhang, Zhimin and Zhou, Zhi},
	journal={SIAM Journal on Scientific Computing},
	volume={42},
	number={2},
	pages={A849--A877},
	year={2020},
	publisher={SIAM}
}

@article{vcertik2013dftatom,
	title={dftatom: A robust and general Schr{\"o}dinger and Dirac solver for atomic structure calculations},
	author={{\v{C}}ert{\'\i}k, Ond{\v{r}}ej and Pask, John E and Vack{\'a}{\v{r}}, Ji{\v{r}}{\'\i}},
	journal={Computer Physics Communications},
	volume={184},
	number={7},
	pages={1777--1791},
	year={2013},
	publisher={Elsevier}
}

@article{wilson2006purgatorio,
	title={Purgatorio—a new implementation of the Inferno algorithm},
	author={Wilson, B and Sonnad, V and Sterne, P and Isaacs, W},
	journal={Journal of Quantitative Spectroscopy and Radiative Transfer},
	volume={99},
	number={1-3},
	pages={658--679},
	year={2006},
	publisher={Elsevier}
}

@article{hohenberg1964inhomogeneous,
	title={Inhomogeneous electron gas},
	author={Hohenberg, Pierre and Kohn, Walter},
	journal={Physical review},
	volume={136},
	number={3B},
	pages={B864},
	year={1964},
	publisher={APS}
}

@article{salvat1991accurate,
	title={Accurate numerical solution of the Schr{\"o}dinger and Dirac wave equations for central fields},
	author={Salvat, Francesc and Mayol, Ricardo},
	journal={Computer physics communications},
	volume={62},
	number={1},
	pages={65--79},
	year={1991},
	publisher={Elsevier}
}

@article{press2007numerical,
	title={Numerical recipes 3rd edition},
	author={Press, William H and Teukolsky, Saul A and Vetterling, William T and Flannery, Brian P},
	journal={Cambridge: New York},
	year={2007}
}

@article{jonsson2007grasp2k,
	title={The grasp2K relativistic atomic structure package},
	author={J{\"o}nsson, Per and He, X and Fischer, C Froese and Grant, IP},
	journal={Computer Physics Communications},
	volume={177},
	number={7},
	pages={597--622},
	year={2007},
	publisher={Elsevier}
}

@article{pask2005finite,
	title={Finite element methods in ab initio electronic structure calculations},
	author={Pask, JE and Sterne, PA},
	journal={Modelling and Simulation in Materials Science and Engineering},
	volume={13},
	number={3},
	pages={R71},
	year={2005},
	publisher={IOP Publishing}
}

@article{grant2009b,
	title={B-spline methods for radial Dirac equations},
	author={Grant, IP},
	journal={Journal of Physics B: Atomic, Molecular and Optical Physics},
	volume={42},
	number={5},
	pages={055002},
	year={2009},
	publisher={IOP Publishing}
}

@article{fischer2009b,
	title={A B-spline Galerkin method for the Dirac equation},
	author={Fischer, Charlotte Froese and Zatsarinny, Oleg},
	journal={Computer Physics Communications},
	volume={180},
	number={6},
	pages={879--886},
	year={2009},
	publisher={Elsevier}
}

\end{document}